\documentclass[12pt]{article}
\usepackage{epsfig}
\usepackage{algorithm,algorithmic}
\usepackage{amsmath}
\usepackage{amssymb}
\usepackage{url}

\usepackage{relsize}

\usepackage{enumitem}

\begin{document}

\def\proof{\noindent {\sl Proof.\ \ }}
\def\endproof{\hfill$\square$
\medskip}

\def\tr{{\rm tr}}
\def\arg{{\rm arg}}
\title{\vspace{-17mm}
Characterization and Computation of Matrices of Maximal Trace over Rotations}
\author{Javier Bernal$^1$, Jim Lawrence$^{1,2}$\\
$^1${\small \sl National Institute of Standards and Technology,} \\
{\small \sl Gaithersburg, MD 20899, USA} \\
$^2${\small \sl George Mason University,} \\
{\small \sl 4400 University Dr, Fairfax, VA 22030, USA} \\
{\tt\small $\{$javier.bernal,james.lawrence$\}$@nist.gov \ \ lawrence@gmu.edu}}
\date{\ }
\maketitle
\begin{abstract}
The constrained orthogonal Procrustes problem is the least-squares problem that calls for
a rotation matrix that optimally aligns two corresponding sets of points in $d-$dimensional
Euclidean space. This problem generalizes to the so-called Wahba's problem which is the
same problem with nonnegative weights. Given a $d\times d$ matrix $M$, solutions to these problems
are intimately related to the problem of finding a $d\times d$ rotation matrix $U$ that maximizes
the trace of $UM$, i.e., that makes $UM$ a matrix of maximal trace over rotation matrices, and
it is well known this can be achieved with a method based on the computation of the singular
value decomposition (SVD) of~$M$. As the main goal of this paper, we characterize $d\times d$
matrices of maximal trace over rotation matrices in terms of their eigenvalues, and for
$d = 2, 3$, we show how this characterization can be used to determine whether a matrix is
of maximal trace over rotation matrices. Finally, although depending only slightly on the
characterization, as a secondary goal of the paper, for $d = 2, 3$, we identify alternative
ways, other than the SVD, of obtaining solutions to the aforementioned problems.
%
\\[0.2cm]
\end{abstract}

\tableofcontents

\label{first}

\section{\normalsize Introduction}
Suppose $P=\{p_1,\ldots,p_n\}$ and $Q=\{q_1,\ldots,q_n\}$ are each sets of $n$ points
in~$\mathbb{R}^d$. With \mbox{$\|\cdot\|$} denoting the $d-$dimensional Euclidean norm,
in the {\em constrained orthogonal Procrustes problem} \cite{kabsch1,kabsch2,umeyama},
a $d\times d$ orthogonal matrix $U$ is found that minimizes
$\Delta(P,Q,U)=\sum_{i=1}^n \|Uq_i-p_i\|^2$, $U$ constrained to be a rotation matrix,
i.e., an orthogonal matrix of determinant~one.
This problem generalizes to the so-called {\em Wahba's problem}
\cite{wahba,markley} which is that of finding a $d\times d$ rotation matrix $U$ that minimizes
$\Delta(P,Q,W,U)=\sum_{i=1}^n\,w_i\,\|Uq_i-p_i\|^2$, where $W=\{w_1,\ldots,w_n\}$ is a set of
$n$ nonnegative weights. Solutions to these problems are of importance, notably in the field
of functional and shape data analysis \cite{srivastava,dogan}, where, in particular, the shapes
of two curves are compared, in part by optimally rotating one curve to match the other.
In \cite{bernal}, for the same purpose, Wahba's problem
occurs with the additional constraint $\sum_{i=1}^n\,w_i=1$ due to the use of the trapezoidal
rule during the discretization of integrals. Given a $d\times d$ matrix $M$,
it is well known that solutions to these problems are intimately related to the problem of
finding among all $d\times d$ rotation matrices~$U$, one that maximizes the trace of~$UM$,
and that the maximization can be achieved with a method called the Kabsch-Umeyama algorithm
(loosely referred to as ``the SVD method'' in what follows), based on the computation of the
{\em singular value decomposition} (SVD) of~$M$ \cite{kabsch1,kabsch2,umeyama,lawrence,
markley,bernal}. In this paper, we analyze matrices of maximal trace over rotation matrices:
A $d\times d$ matrix $M$ is of maximal trace over rotation matrices if given any $d\times d$
rotation matrix~$U$, the trace of $UM$ does not exceed that of~$M$. As a result, we identify
a characterization of these matrices: A $d\times d$ matrix is of maximal trace over rotation
matrices if and only if it is symmetric and has at most one negative eigenvalue, which, if it
exists, is no larger in absolute value than the other eigenvalues of the matrix.
Establishing this characterization is the main goal of this paper, and for $d=2,3$, we show
how it can be used to determine whether a matrix is of maximal trace over
rotation matrices. Finally, although depending only slightly on the characterization, as a
secondary goal of the paper, for $d=2,3$, we identify alternative ways, other than the SVD,
of obtaining solutions to the aforementioned problems. Accordingly, for $d=2$, we identify an
alternative way that does not involve the SVD, and for $d=3$, one that without it, for matrices
of randomly generated entries, is successful in our experiments close to one hundred percent
of the time, using the SVD only when it is not.

\par In Section~2 we reformulate the constrained orthogonal Procrustes problem and Wahba's problem
in terms of the trace of matrices, and verify the well-known fact that given one such problem,
there is a $d\times d$ matrix $M$ associated with it such that a $d\times d$ rotation matrix $U$
is a solution to it if and only if $UM$ is of maximal trace over rotation matrices.
In Section~3 we identify the characterization of $d\times d$ matrices of maximal
trace over rotation matrices and show that for $d=2,3$, it can be used to determine whether
a matrix is of maximal trace over rotation matrices. Once the main goal of the paper is
established, i.e., the characterization has been identified, most of the rest of Section~3 and
for that matter the rest of the paper, is for accomplishing the secondary goal of the paper,
i.e., identifying, for $d=2,3$, alternative ways, other than the SVD, of obtaining solutions
to the aforementioned problems.
For this purpose, in Section~4 we present alternative solutions expressed in closed form, that
do not involve the~SVD method, to the two-dimensional constrained orthogonal Procrustes problem
and Wahba's problem, and indeed, given a $2\times 2$ matrix $M$, to the problem of finding a
$2\times 2$ rotation matrix $U$ such that $UM$ is of maximal trace over rotation matrices.
Using results from the latter part of Section~3, in~Section~5, given a $3\times 3$ symmetric
matrix $M$, we present an alternative solution that does not involve the~SVD to the problem
of finding a rotation matrix $U$ such that $UM$ is of maximal trace over rotation matrices.
This alternative solution is based on a trigonometric identity that can still be used
if the matrix $M$ is not symmetric, to produce the usual orthogonal matrices necessary to
carry out the SVD method. Finally, in Section~6, we reconsider the situation in which the
$3\times 3$ matrix $M$ is not symmetric, and as an alternative to the SVD method present a
procedure that uses the so-called Cayley transform in conjunction with Newton's method to
find a $3\times 3$ rotation matrix $U$ so that $UM$ is symmetric, possibly of maximal trace
over rotation matrices. If the resulting $UM$ is not of maximal trace over rotation matrices,
using the fact that it is symmetric, another $3\times 3$ rotation matrix $R$ can then be
computed (without the SVD) as described in Section~5 so that $RUM$ is of maximal trace over
rotation~matrices. Of course, if Newton's method fails in the procedure, the SVD method is
still used as described above. We then note, still in Section~6, that all of the above about
the three-dimensional case, including carrying out the SVD method as described above, has been
successfully implemented in Fortran, and that without the SVD, for randomly generated matrices,
the Fortran code is successful in our experiments close to one hundred per cent of the time,
using the SVD only when it is not. We then also note that the Fortran code is faster than
Matlab\footnote{The identification of any commercial product or trade name does not imply
endorsement or recommendation by the National Institute of Standards and Technology.} code
using Matlab's SVD command, and provide links to all codes (Fortran and Matlab) at the end
of the section.
\section{\normalsize Reformulation of Problems as Maximizations of the Trace of Matrices
Over Rotations}
With $P=\{p_1,\ldots,p_n\}$, $Q=\{q_1,\ldots,q_n\}$, each a set of $n$ points
in~$\mathbb{R}^d$, and $W=\{w_1,\ldots,w_n\}$, a set of $n$ nonnegative real numbers (weights),
we now think of $P$ and $Q$ as $d\times n$ matrices having the $p_i$'s and $q_i$'s
as columns so that $P=\big(p_1 \ldots p_n\big)$, $Q=\big(q_1 \ldots q_n\big)$, and of $W$ as an
$n\times n$ diagonal matrix with $w_1,\ldots,w_n$ as the elements of the diagonal,
in that order running from the upper left to the lower right of~$W$ so that
%
\[W= \begin{pmatrix} w_1 & 0 & \ldots & \ldots & 0 \cr
                     0 & w_2 & 0 & \ldots & 0 \cr
                     \ldots & \ldots & \ldots & \ldots & \ldots \cr
                     0 & \ldots & 0 & w_{n-1} & 0 \cr
                     0 & \ldots & \ldots & 0 & w_n \cr \end{pmatrix}.\]
\par
Since Wahba's problem becomes the constrained orthogonal Procrustes problem if the weights are all
set to~one, we focus our attention on Wahba's problem and thus wish to find a $d\times d$ rotation
matrix $U$ that minimizes $$\Delta(P,Q,W,U)=\sum_{i=1}^n\,w_i\,\|Uq_i-p_i\|^2.$$ With this purpose
in mind, we rewrite $\sum_{i=1}^n\,w_i\,\|Uq_i-p_i\|^2$ as follows, where given a matrix $R$,
$\tr(R)$ stands for the trace of~$R$
\begin{eqnarray*}
&&\hspace*{-.325in} \sum_{i=1}^n w_i\,||Uq_i -p_i||^2 = \sum_{i=1}^n w_i(Uq_i - p_i)^T (Uq_i - p_i)\\
&=& \tr\bigl(W(UQ - P)^T(UQ-P)\bigr)\\
&=& \tr\bigl(W(Q^TU^T - P^T) (UQ - P)\bigr)\\
&=&\tr\big(W(Q^TQ + P^TP - Q^TU^TP - P^T UQ)\big)\\
&=& \tr(WQ^TQ) + \tr(WP^TP)- \tr(WQ^TU^TP) - \tr(WP^TUQ)\\
&=& \tr(WQ^TQ) + \tr(WP^TP)- \tr(P^TUQW^T) - \tr(WP^TUQ)\\
&=& \tr(WQ^TQ) + \tr(WP^TP)- \tr(W^TP^TUQ) - \tr(WP^TUQ)\\
&=& \tr(WQ^TQ) + \tr(WP^TP) - 2\tr(WP^TUQ)\\
&=& \tr(WQ^TQ) + \tr(WP^TP) - 2\tr(UQWP^T)
\end{eqnarray*}
where a couple of times we have used the fact that for positive integers $k$, $l$,
if $A$ is a $k\times l$ matrix and $B$ is an $l\times k$ matrix, then $\tr(AB)=\tr(BA)$.
Since only the third term in the last line above depends on $U$, denoting the $d\times d$
matrix $QWP^T$ by~$M$, it suffices to find a $d\times d$ rotation matrix $U$ that
maximizes~$\tr(UM)$, and it is well known that one such $U$ can be computed with a
method based on the singular value decomposition (SVD) of~$M$
\cite{kabsch1,kabsch2,umeyama,lawrence, markley,bernal}. This method is called the
Kabsch-Umeyama algorithm~\cite{kabsch1,kabsch2,umeyama} (loosely referred to as
``the SVD method'' throughout this paper), which we outline for the sake of completeness
(see~Algorithm Kabsch-Umeyama below, where $\mathrm{diag}\{s_1,\ldots,s_d\}$
is the $d\times d$ diagonal matrix with numbers $s_1,\ldots,s_d$ as the elements
of the diagonal, in that order running from the upper left to the lower right of the matrix,
and see~\cite{lawrence} for its justification in a purely algebraic manner through the
exclusive use of simple concepts from linear algebra).
The SVD~\cite{lay} of~$M$ is a representation of the form $M=VSR^T$, where $V$ and $R$
are~$d\times d$ orthogonal matrices and $S$ is a $d\times d$ diagonal matrix with the singular
values of $M$, which are nonnegative real numbers, appearing in the diagonal of $S$ in descending
order, from  the upper left to the lower right of~$S$. Finally, note that any real matrix,
not necessarily square, has an~SVD , not necessarily unique~\cite{lay}.
\begin{algorithmic}
\STATE \noindent\rule{13cm}{0.4pt}
\STATE {\bf Algorithm Kabsch-Umeyama}
\STATE \noindent\rule[.1in]{13cm}{0.4pt}
\STATE Compute $d\times d$ matrix $M=QWP^T$.
\STATE Compute SVD of $M$, i.e., identify $d\times d$ matrices $V$, $S$, $R$,
\STATE with $M = V S R^T$ in the SVD sense.
\STATE Set $s_1= \ldots = s_{d-1}=1$.
\STATE If $\det(VR) > 0$, then set $s_d=1$, else set $s_d=-1$.
\STATE Set $\tilde{S} = \mathrm{diag}\{s_1,\ldots,s_d\}$.
\STATE Return $d\times d$ rotation matrix $U = R \tilde{S} V^T$.
\STATE \noindent\rule{13cm}{0.4pt}
\end{algorithmic}
\smallskip
\section{\normalsize Characterization of Matrices of Maximal Trace Over Rotations}
In what follows, given a real $d\times d$ matrix $M$, we say that $M$ {\em is of maximal trace over
rotation (orthogonal) matrices} if for any $d\times d$ rotation (orthogonal) matrix $U$,
it must be that~$\tr(M)\geq\tr(UM)$.\,\medskip\\
{\bf Proposition 1:} Let $M$ be a $d\times d$ matrix. If one of the following occurs,
then the other two occur as well.
\begin{enumerate}[nolistsep]
 \item[i.] $\tr(M)\geq\tr(UM)$ for any $d\times d$ rotation (orthogonal) matrix~$U$.
 \item[ii.] $\tr(M)\geq\tr(MU)$ for any $d\times d$ rotation (orthogonal) matrix~$U$.
 \item[iii.] $\tr(M)\geq\tr(UMV)$ for any $d\times d$ rotation (orthogonal) matrices~$U, V$.
\end{enumerate}
\,\smallskip
{\bf Proof:} With $I$ as the $d\times d$ identity matrix, we have\\
\hspace*{.205in}i $\Rightarrow$\ \ ii: $\tr(MU)=\tr(UMUU^T)=\tr(UM)\leq\tr(M)$.\\
\hspace*{.165in}ii $\Rightarrow$\;\,iii: $\tr(UMV)=\tr(U^TUMVU)=\tr(MVU)\leq\tr(M)$.\\
\hspace*{.135in}iii $\Rightarrow$\ \ \;i: $\tr(UM)=\tr(UMI)\leq\tr(M)$.\\
\smallskip \endproof
\,\medskip\\
{\bf Proposition 2:} Let $A$ be a $d\times d$ matrix. If $A$ is of maximal trace over
rotation (orthogonal) matrices, then $A$ is a symmetric matrix.\,\medskip\\
{\bf Proof:}
It suffices to prove the proposition only for rotation matrices.\\
Let $a_{ij}$, $i,j=1,\ldots,d$, be the entries in~$A$, and assuming $A$ is not
a symmetric matrix, suppose $k,l$ are such that $k>l$ and $a_{lk}\not=a_{kl}$.\\
Given an angle $\theta$, $0\leq\theta<2\pi$, a $d\times d$ so-called Givens rotation
$G(k,l,\theta)$ can be defined with entries $g_{ij}$, $i,j=1,\dots,d$, among which, the nonzero entries
are given by
\begin{eqnarray*}
g_{mm}&=&1,\ m=1,\ldots,d,\ m\not=l,\ m\not=k\\
g_{ll}&=&g_{kk}=\cos\theta\\
g_{lk}&=&-g_{kl}= -\sin\theta.
\end{eqnarray*}
That $G(k,l,\theta)$ is a rotation matrix (an orthogonal matrix of determinant equal to one) has long
been established and is actually easy to verify.\\
Let $a=a_{ll}+a_{kk}$, $b=a_{lk}-a_{kl}$, $c=\sqrt{a^2+b^2}$.
Clearly $b\not=0$ so that~$c\not=0$ and~$c>a$.
For our purposes we choose $\theta$ so that $\cos\theta=a/c$ and $\sin\theta=b/c$.
Thus, for this $\theta$, $G(k,l,\theta)$ is such that $g_{ll}=g_{kk}=a/c$, and $g_{lk}=-g_{kl}=-b/c$.\\
We show $\tr(G(k,l,\theta)A)>\tr(A)$ which contradicts that $A$ is of maximal trace over rotation
matrices.\\
Let $v_{ij}$, $i,j=1,\ldots,d$, be the entries in~$G(k,l,\theta)A$.
We show $\sum_{m=1}^d v_{mm} > \sum_{m=1}^d a_{mm}$.
Clearly $v_{mm}=a_{mm}$, $m=1,\ldots,d$, $m\not=l$, $m\not=k$, thus it suffices to show
$v_{ll}+v_{kk}>a_{ll}+a_{kk}$, and we know $a_{ll}+a_{kk}=a$.\\
Also clearly $v_{ll} = g_{ll}a_{ll}+g_{lk}a_{kl}$ and $v_{kk}= g_{kl}a_{lk}+g_{kk}a_{kk}$, so that
\begin{eqnarray*}
v_{ll}+v_{kk}&=&(a/c)a_{ll}+(-b/c)a_{kl}+(b/c)a_{lk}+(a/c)a_{kk}\\
&=&(a/c)(a_{ll}+a_{kk})+(b/c)(a_{lk}-a_{kl})=(a/c)a+(b/c)b\\
&=&(a^2+b^2)/c=c^2/c=c>a.
\end{eqnarray*}
\endproof
\,\medskip
\par The following useful proposition was proven in \cite{lawrence}. For the sake of
completeness we present the proof here. Here
$\mathrm{diag}\{\sigma_1,\ldots,\sigma_d\}$ is the $d\times d$ diagonal matrix with the numbers
$\sigma_1,\ldots,\sigma_d$ as the elements of the diagonal, in that order running from the upper
left to the lower right of the~matrix, and det$(W)$ is the determinant of the $d\times d$ matrix~$W$.
\,\medskip\\
{\bf Proposition 3:} If $D=\mathrm{diag}\{\sigma_1,\ldots,\sigma_d\}$, $\sigma_j\geq 0$,
$j=1,\ldots,d$, and $W$ is a $d\times d$ orthogonal matrix, then
\begin{enumerate}[nolistsep]
\item [$1)$] tr$(WD)\leq\sum_{j=1}^d \sigma_j$.
\item [$2)$] If $B$ is a $d\times d$ orthogonal matrix, $S=B^TDB$, then $\tr(WS)\leq \tr(S)$.
\item [$3)$] If det$(W)=-1$, $\sigma_d\leq \sigma_j$, $j=1,\ldots,d-1$, then
$\tr(WD)\leq\sum_{j=1}^{d-1}\sigma_j-\sigma_d$.
\end{enumerate}
\,\smallskip
{\bf Proof:} Since $W$ is orthogonal and if $W_{kj}$, $k,j=1,\ldots,d$, are the
entries of $W$, then, in particular, $W_{jj}\leq 1$, $j=1,\ldots,d$, so that\\
$\mathrm{tr}(WD)=\sum_{j=1}^d W_{jj}\sigma_j\leq \sum_{j=1}^d \sigma_j$, and therefore statement~$1)$
holds.\\
Accordingly, assuming $B$ is a $d\times d$ orthogonal matrix, since $BWB^T$ is also
orthogonal, it follows from~$1)$ that\\
$\mathrm{tr}(WS)=\mathrm{tr}(WB^TDB)=\mathrm{tr}(BWB^TD)\leq\sum_{j=1}^d\sigma_j=\mathrm{tr}(D)
=\mathrm{tr}(S)$, and therefore~$2)$ holds.\\
If det$(W)=-1$, we show next that a $d\times d$ orthogonal matrix $B$ can be
identified so that with $\bar{W}=B^TWB$, then
$\bar{W}= \left( \begin{smallmatrix}
W_0 & O\\ O^T & -1\\
\end{smallmatrix} \right)$,
in which $W_0$ is interpreted as the upper leftmost $d-1\times d-1$ entries of $\bar{W}$ and as a
$d-1\times d-1$ matrix as well, and $O$ is interpreted as a vertical column or vector of $d-1$ zeroes.\\
With $I$ as the $d\times d$ identity matrix, then det$(W)=-1$ implies
$\mathrm{det}(W+I)=-\mathrm{det}(W)\mathrm{det}(W+I)=-\mathrm{det}(W^T)\mathrm{det}(W+I)=
-\mathrm{det}(I+W^T)=-\mathrm{det}(I+W)$ which implies det$(W+I)=0$
so that $x\not=0$ exists in $\mathbb{R}^d$ with $Wx=-x$. It also follows then that
$W^TWx=W^T(-x)$ which gives $x=-W^Tx$ so that $W^Tx=-x$ as well.\\
Letting $b_d=x$, vectors $b_1,\ldots,b_{d-1}$ can be obtained so that $b_1,\ldots,b_d$ form
a basis of~$\mathbb{R}^d$, and by the Gram-Schmidt process starting with $b_d$, we may
assume $b_1,\ldots,b_d$ form an orthonormal basis of $\mathbb{R}^d$ with $Wb_d=W^Tb_d=-b_d$.
Letting $B=(b_1,\ldots,b_d)$, interpreted as a $d\times d$ matrix with columns $b_1,\ldots,b_d$,
in that order, it then follows that $B$ is orthogonal, and with $\bar{W}=B^TWB$ and
$W_0$, $O$ as previously described, noting
$B^TWb_d=B^T(-b_d)= \left( \begin{smallmatrix} O\\ -1\\ \end{smallmatrix} \right)$ and
$b_d^TWB=(W^Tb_d)^TB=(-b_d)^TB=(O^T \, -1)$,
then $\bar{W}= \left( \begin{smallmatrix}
W_0 & O\\ O^T & -1\\
\end{smallmatrix} \right)$.
Note $\bar{W}$ is orthogonal and therefore so is the $d-1\times d-1$ matrix~$W_0$.\\
Let $S=B^TDB$ and write
$S = \left( \begin{smallmatrix}
S_0 & a\\ b^T & \gamma\\
\end{smallmatrix} \right)$,
in which $S_0$ is interpreted as the upper leftmost $d-1\times d-1$ entries of $S$ and as a
$d-1\times d-1$ matrix as well, $a$ and $b$ are interpreted as vertical columns or vectors of
$d-1$ entries, and $\gamma$ as a scalar.\\
Note $\mathrm{tr}(WD)=\mathrm{tr}(B^TWDB)=\mathrm{tr}(B^TWBB^TDB)=\mathrm{tr}(\bar{W}S)$,
so that $\bar{W}S=$
$\left( \begin{smallmatrix}
W_0 & O\\ O^T & -1\\
\end{smallmatrix} \right)$
$\left( \begin{smallmatrix}
S_0 & a\\ b^T & \gamma\\
\end{smallmatrix} \right)=$
$\left( \begin{smallmatrix}
W_0S_0 & W_0a\\ -b^T & -\gamma\\
\end{smallmatrix} \right)$
gives $\mathrm{tr}(WD)=\mathrm{tr}(W_0S_0)-\gamma$.\\
We show $\mathrm{tr}(W_0S_0)\leq\mathrm{tr}(S_0)$. For this purpose let
$\hat{W}= \left( \begin{smallmatrix}
W_0 & O\\ O^T & 1\\
\end{smallmatrix} \right)$,
$W_0$ and $O$ as above. Since $W_0$ is orthogonal, then clearly $\hat{W}$ is a $d\times d$
orthogonal matrix, and by~$2)$, $\mathrm{tr}(\hat{W}S)\leq \mathrm{tr}(S)$
so that $\hat{W}S=$
$\left( \begin{smallmatrix}
W_0 & O\\ O^T & 1\\
\end{smallmatrix} \right)$
$\left( \begin{smallmatrix}
S_0 & a\\ b^T & \gamma\\
\end{smallmatrix} \right)=$
$\left( \begin{smallmatrix}
W_0S_0 & W_0a\\ b^T & \gamma\\
\end{smallmatrix} \right)$
gives $\mathrm{tr}(W_0S_0)+\gamma=\mathrm{tr}(\hat{W}S)\leq\mathrm{tr}(S)=\mathrm{tr}(S_0)+\gamma$.
Thus, $\mathrm{tr}(W_0S_0)\leq\mathrm{tr}(S_0)$.\\
Note $\mathrm{tr}(S_0)+\gamma = \mathrm{tr}(S)= \mathrm{tr}(D)$, and if $B_{kj}$,
$k,j=1,\ldots,d$ are the entries of~$B$, then $\gamma = \sum_{k=1}^d B_{kd}^2\sigma_k$,
a convex combination of the $\sigma_k$'s, so that $\gamma\geq\sigma_d$.
It then follows that\\ $\mathrm{tr}(WD)=\mathrm{tr}(W_0S_0)-\gamma\leq\mathrm{tr}(S_0)-\gamma=
\mathrm{tr}(D)-\gamma-\gamma\leq\sum_{j=1}^{d-1}\sigma_j-\sigma_d$, and therefore~$3)$ holds.
\endproof
\,\smallskip
\par
Conclusion~$3)$ of Proposition 3 above can be improved as follows.\,\medskip\\
{\bf Proposition 4:} Given $D=\mathrm{diag}\{\sigma_1,\ldots,\sigma_d\}$, $\sigma_j\geq 0$,
$j=1,\ldots,d$, let $k = \arg\,\min_j\{\sigma_j,j=1,\ldots,d\}$. If $W$ is a $d\times d$ orthogonal
matrix with det$(W)=-1$, then $\tr(WD)\leq\sum_{j=1,j\not=k}^{d}\sigma_j-\sigma_k$.\,\medskip\\
{\bf Proof:} Assume $k\not=d$, as otherwise the result follows from~$3)$ above.\\
Let $G$ be the $d\times d$ orthogonal matrix with entries $g_{il}$, $i,l=1\ldots,d$,
among which, the nonzero entries are given by
\begin{eqnarray*}
g_{mm} = 1,\ m=1,\ldots,d,\ m\not=k,\ m\not=d,\ \ g_{kd}=g_{dk}=1.
\end{eqnarray*}
Note $G^{-1}=G^T=G$. Letting $\hat{W}=GWG$, $\hat{D}=GDG$, then det$(\hat{W})=-1$ and
$\tr(\hat{W}\hat{D})=\tr(G^TWGG^TDG)=\tr(WD)$.\\
Note $\hat{W}$ is the result of switching rows $k$ and $d$ of $W$ and then switching columns $k$
and $d$ of the resulting matrix. The same applies to~$\hat{D}$ with respect to $D$ so that $\hat{D}$
is a diagonal matrix whose diagonal is still the diagonal of $D$ but with $\sigma_k$ and $\sigma_d$
trading places in~it. It follows then by~$3)$ of Proposition~3 that
$\tr(WD)=\tr(\hat{W}\hat{D})= \leq\sum_{j=1,j\not=k}^{d}\sigma_j-\sigma_k$.
\endproof
\,\medskip\par
In the following two propositions and corollary, matrices of maximal trace over rotation (orthogonal)
matrices are characterized.
\,\medskip\\
{\bf Proposition 5:} If $A$ is a $d\times d$ symmetric matrix, then
\begin{enumerate}[nolistsep]
\item[$1)$] If $A$ is positive semidefinite (this is equivalent to each eigenvalue
of $A$ being nonnegative), then det$(A)\geq 0$ and $A$ is of maximal trace over orthogonal
matrices, and therefore over rotation matrices.
\item[$2)$] If $A$ has exactly one negative eigenvalue, the absolute value of this eigenvalue
being at most as large as any of the other eigenvalues, then det$(A)<0$ and $A$ is of maximal trace
over rotation matrices.
\end{enumerate}
\,\smallskip
{\bf Proof:} Since $A$ is a symmetric matrix, there are $d\times d$ matrices $V$ and $D$, real
numbers $\alpha_j$, $j=1,\ldots,d$, $V$ orthogonal, $D=\mathrm{diag}\{\alpha_1,\ldots,\alpha_d\}$,
with $A=V^TDV$ so~that $\tr(A)=\tr(D)=\sum_{j=1}^d\alpha_j$,
the~$\alpha_j$'s the eigenvalues of~$A$.\medskip\\
If $A$ is positive semidefinite, then $\alpha_j\geq 0$, $j=1,\ldots,d$, and det$(A)\geq 0$.\\
Let $W$ be a $d\times d$ orthogonal matrix. Then by $1)$ of Proposition~3\\
$\tr(WA)=\tr(WV^TDV)=\tr(VWV^TD)\leq\sum_{j=1}^d\alpha_j=\tr(A)$\\ and therefore
statement~$1)$ holds.\medskip\\
If $A$ has exactly one negative eigenvalue, the absolute value of this eigenvalue being
at most as large as any of the other eigenvalues, let $k$, $1\leq k\leq d$, be such that
$\alpha_k<0$, $|\alpha_k|\leq\alpha_j$, $j=1,\ldots,d$, $j\not=k$. Clearly det$(A)<0$.\\
Let $\sigma_j=\alpha_j$, $j=1\ldots,d$, $j\not=k$, $\sigma_k=-\alpha_k$, and
$\hat{D}=\mathrm{diag}\{\sigma_1,\ldots,\sigma_d\}$.\\
Let $G$ be the orthogonal matrix with entries $g_{il}$, $i,l=1,\ldots,d$, among which,
the nonzero entries are given by
\begin{eqnarray*}
g_{mm}=1,\ m=1,\ldots,d,\ m\not=k,\ \ g_{kk}=-1.
\end{eqnarray*}
Note det$(G)=-1$, $G^{-1}=G^T=G$, $\hat{D}=GD$.\\
Let $U$ be a $d\times d$ rotation matrix. Letting $W=VUV^TG$, then det$(W)=-1$.\\
By Proposition 4, then\\
$\tr(UA)=\tr(UV^TDV)=\tr(VUV^TD)=\tr(VUV^TGGD)=\tr(W\hat{D})\\
\hspace*{.475in}\leq\sum_{j=1,j\not=k}^d\sigma_j-\sigma_k=\sum_{j=1}^d\alpha_j=\tr(A)$ as
$-\sigma_k=\alpha_k$\\
and therefore statement~$2)$ holds.
\endproof
\,\medskip\\
{\bf Proposition 6:}
If $A$ is a $d\times d$ matrix of maximal trace over orthogonal matrices,
then $A$ is a symmetric matrix and as such it is positive semidefinite.\\
On the other hand, if $A$ is a $d\times d$ matrix of maximal trace over
rotation matrices, then $A$ is a symmetric matrix and
\begin{enumerate}[nolistsep]
\item[$1)$] If det$(A)=0$, then $A$ is positive semidefinite (this is equivalent to each eigenvalue
of $A$ being nonnegative).
\item[$2)$] If det$(A)>0$, then $A$ is positive definite (this is equivalent to each eigenvalue
of $A$ being positive).
\item[$3)$] If det$(A)<0$, then $A$ has exactly one negative eigenvalue, and the absolute value of
this eigenvalue is at most as large as any of the other eigenvalues.
\end{enumerate}
\,\smallskip
{\bf Proof:} That $A$ is symmetric for all cases follows from Proposition~2.
Accordingly, there are $d\times d$ matrices $V$ and $D$, real numbers $\alpha_j$,
$j=1,\ldots,d$, $V$ orthogonal, $D=\mathrm{diag}\{\alpha_1,\ldots,\alpha_d\}$, with
$A=V^TDV$ so~that $\tr(A)=\tr(D)=\sum_{j=1}^d\alpha_j$, the~$\alpha_j$'s the eigenvalues
of~$A$.\medskip\\
Assume $A$ is of maximal trace over orthogonal matrices and $A$ is not positive semidefinite.\\
Then $A$ must have a negative eigenvalue.\\
Accordingly, let $k$, $1\leq k \leq d$, be such that $\alpha_k<0$, and
let $\sigma_j=\alpha_j$, $j=1,\ldots,d$, $j\not=k$, $\sigma_k=-\alpha_k$,
$\hat{D}=\mathrm{diag}\{\sigma_1,\ldots,\sigma_d\}$.\\
Let $G$ be the orthogonal matrix with entries $g_{ih}$, $i,h=1,\ldots,d$, among which,
the nonzero entries are given by
\begin{eqnarray*}
g_{mm}=1,\ m=1,\ldots,d,\ m\not=k,\ \ g_{kk}=-1.
\end{eqnarray*}
Note $GD=\hat{D}$ so that letting $U=V^TGV$, then $U$ is orthogonal and\\
$\tr(UA)=\tr(V^TGVV^TDV)=\tr(GD)=\tr(\hat{D})=\sum_{j=1}^d\sigma_j
>\sum_{j=1}^d\alpha_j=\tr(A)$ as $\sigma_k>\alpha_k$ which contradicts that
$A$ is of maximal trace over orthogonal matrices. Thus, it must be that $A$ is positive
semidefinite and therefore the first part of the proposition holds.\medskip\\
Assume now $A$ is of maximal trace over rotation matrices.\\
Before proceeding with the rest of the proof, we define some matrices and numbers that are
used repeatedly throughout the proof in the same manner, and make some observations about them.\\
Accordingly, let $k,l$, $k\not=l$, $1\leq k,l \leq d$, be given, and
let $\sigma_j=\alpha_j$, $j=1,\ldots,d$, $j\not=k$, $j\not=l$,
$\sigma_k=-\alpha_k$, $\sigma_l=-\alpha_l$,
$\hat{D}=\mathrm{diag}\{\sigma_1,\ldots,\sigma_d\}$.\\
Let $G$ be the orthogonal matrix with entries $g_{ih}$, $i,h=1,\ldots,d$, among which,
the nonzero entries are given by
\begin{eqnarray*}
g_{mm}=1,\ m=1,\ldots,d,\ m\not=k,\ m\not=l,\ \ g_{kk}=g_{ll}=-1.
\end{eqnarray*}
Note det$(G)=1$ and $GD=\hat{D}$ so that letting $U=V^TGV$, then det$(U)=1$\\ and
$\tr(UA)=\tr(V^TGVV^TDV)=\tr(GD)=\tr(\hat{D})=\sum_{j=1}^d\sigma_j$.\medskip\\
If det$(A)=0$, assume $A$ is not positive semidefinite.
Then $A$ must have an eigenvalue equal to zero and a negative eigenvalue.
Accordingly, let $k,l$, $k\not=l$, $1\leq k,l \leq d$, be such that $\alpha_k=0$, $\alpha_l<0$.\\
With $U$ and $\sigma_j$, $j=1,\dots,d$, as defined above, note $\sigma_k=\alpha_k=0$,
$\sigma_l>\alpha_l$, so that $\tr(UA)=\sum_{j=1}^d\sigma_j>\sum_{j=1}^d\alpha_j=\tr(A)$
which contradicts that $A$ is of maximal trace over rotation matrices. Thus, it must be that $A$
is positive semidefinite and therefore statement~$1)$ holds.\medskip\\
If det$(A)>0$, assume $A$ is not positive definite.
Then $A$ must have an even number of negative eigenvalues.
Accordingly, let $k,l$, $k\not= l$, $1\leq k,l \leq d$, be such that $\alpha_k<0$, $\alpha_l<0$.\\
With $U$ and $\sigma_j$, $j=1,\dots,d$, as defined above, note $\sigma_k>\alpha_k$,
$\sigma_l>\alpha_l$, so that $\tr(UA)=\sum_{j=1}^d\sigma_j>\sum_{j=1}^d\alpha_j=\tr(A)$
which contradicts that $A$ is of maximal trace over rotation matrices. Thus, it must be that $A$
is positive definite and therefore statement~$2)$ holds.\medskip\\
If det$(A)<0$, then $A$ has at least one negative eigenvalue.
Assume first $A$ has more than one negative eigenvalue.
Accordingly, let $k,l$, $k\not=l$, $1\leq k,l \leq d$, be such that $\alpha_k<0$, $\alpha_l<0$.\\
With $U$ and $\sigma_j$, $j=1,\dots,d$, as defined above, note $\sigma_k>\alpha_k$,
$\sigma_l>\alpha_l$, so that $\tr(UA)=\sum_{j=1}^d\sigma_j>\sum_{j=1}^d\alpha_j=\tr(A)$
which contradicts that $A$ is of maximal trace over rotation matrices. Thus, it must be that $A$
has exactly one negative eigenvalue.\\
Assume now that the absolute value of the only negative eigenvalue of $A$ is larger than
some other (nonnegative) eigenvalue of~$A$.
Accordingly, let $k, l$, $k\not=l$, $1\leq k \leq d$, be such that $\alpha_k<0$ so that
$\alpha_k$ is the only negative eigenvalue of~$A$, and $|\alpha_k|>\alpha_l\geq 0$.\\
With $U$ and $\sigma_j$, $j=1,\dots,d$, as defined above, note
$\sigma_k+\sigma_l >0 >\alpha_k+\alpha_l$,
so that $\tr(UA)=\sum_{j=1}^d\sigma_j>\sum_{j=1}^d\alpha_j=\tr(A)$
which contradicts that $A$ is of maximal trace over rotation matrices. Thus, it must be that
$A$ has exactly one negative eigenvalue, and the absolute value of
this eigenvalue is at most as large as any of the other eigenvalues,
and thus statement~$3)$ holds. 
\endproof
\,\medskip\\
{\bf Corollary 1:}
Let $A$ be a $d\times d$ matrix. $A$ is of maximal trace over orthogonal matrices
if and only if $A$ is symmetric and as such it is positive semidefinite.\\
On the other hand, $A$ is of maximal trace over rotation matrices if and only if $A$ is
symmetric and has at most one negative eigenvalue, which, if it exists, is no larger in
absolute value than the other eigenvalues of~$A$. Consequently, if $A$ is of maximal trace over
rotation (orthogonal) matrices, then the trace of $A$ is nonnegative.\,\medskip\\
{\bf Proof:} The sufficiency and necessity of the two characterizations follow from Proposition~5
and Proposition~6, respectively. The last part follows from the characterizations and the fact
that the trace of any matrix equals the sum of its eigenvalues.
\endproof
\,\medskip\par
We note the characterization above involving orthogonal matrices is well known. See
\cite{horn}, page~432. We also note the characterization above involving rotation matrices
is the main result of this paper.\,\medskip\par
Due to the characterization above involving rotation matrices, Proposition~7 and Proposition~8
that follow, provide, respectively, ways of determining whether a symmetric matrix is of maximal
trace over rotation matrices for~$d=2$ and~$d=3$.
\,\medskip\\
{\bf Proposition 7:} Let $A$ be a $2\times 2$ symmetric matrix. Then the trace of $A$ is
nonnegative if and only if $A$ has at most one negative eigenvalue, which, if it exists,
is no larger in absolute value than the other eigenvalue of~$A$. Thus, the trace of $A$ is
nonnegative if and only if $A$ is of maximal trace over rotation~matrices.\,\medskip\\
{\bf Proof:} Let $\alpha$, $\beta$ be the eigenvalues of $A$. If $\alpha<0$ and $\beta<0$,
then $\alpha+\beta<0$, and if $\alpha<0$, $\beta\geq 0$ and $|\alpha|>\beta$, or
$\beta<0$, $\alpha\geq 0$ and $|\beta|>\alpha$, then
$\alpha+\beta<0$. Also, if $\alpha+\beta<0$, then either $\alpha<0$ and $\beta<0$, or
$\alpha<0$, $\beta\geq 0$ and $|\alpha|>\beta$, or $\beta<0$, $\alpha\geq 0$ and $|\beta|>\alpha$.
It is clear then that $\alpha+\beta$ is nonnegative if and only if at most one of
$\alpha$, $\beta$ is negative, in which case the one that is negative must be at most as large
as the other one in absolute value.\\
The last part of the proposition follows then from Corollary 1.
\endproof
\,\medskip\\
{\bf Proposition 8:} Let $A$ be a $3\times 3$ symmetric matrix and let $S=\tr(A)I-A$, where
$I$ is the $3\times 3$ identity matrix. Then $S$ is positive semidefinite if and only if $A$ has at
most one negative eigenvalue, which, if it exists, is no larger in absolute value that the
other two eigenvalues of~$A$. Thus, $S$ is positive semidefinite if and only if $A$ is of
maximal trace over rotation matrices.\,\medskip\\
{\bf Proof:}. Clearly $S$ is a symmetric matrix. Let $\alpha$, $\beta$, $\gamma$ be the
eigenvalues of $A$. Then the eigenvalues of $S$ are $\alpha+\beta$, $\beta+\gamma$,
$\gamma+\alpha$. We only show $\alpha+\beta$~is. Accordingly, let $w\not=0$ be a point in $\mathbb{R}^3$
such that $Aw=\gamma w$. Then $Sw = (\tr(A)I-A)w=(\alpha+\beta+\gamma)w-\gamma w=(\alpha+\beta)w$.
Thus, $\alpha+\beta$~is.\\
If, say $\alpha<0$ and $\beta<0$, then $\alpha+\beta<0$, and
if, say $\alpha<0$, $\beta\geq 0$ and $|\alpha|>\beta$, or
$\beta<0$, $\alpha\geq 0$ and $|\beta|>\alpha$, then $\alpha+\beta<0$.
Also, if, say $\alpha+\beta<0$, then either $\alpha<0$ and $\beta<0$, or
$\alpha<0$, $\beta\geq 0$ and $|\alpha|>\beta$, or $\beta<0$, $\alpha\geq 0$ and $|\beta|>\alpha$.
It is clear then that $\alpha+\beta$, $\beta+\gamma$, $\gamma+\alpha$ are nonnegative if and
only if at most one of $\alpha$, $\beta$, $\gamma$ is negative, in which case the one that is
negative must be at most as large as the other two in absolute value.\\
The last part of the proposition follows then from Corollary 1.
\endproof
\,\medskip\par
Let $A$ be a $d\times d$ symmetric matrix and let $S=\tr(A)I-A$, where $I$ is the $d\times d$
identity matrix. Proposition~8 shows that for $d=3$ a sufficient and necessary condition for
$A$ to be of maximal trace over rotation matrices is that $S$ be positive semidefinite. The
next proposition shows, in particular, that for any $d$, if $d$ is odd, then a necessary
condition for $A$ to be of maximal trace over rotation matrices is that $S$ be positive
semidefinite.\,\medskip\\
{\bf Proposition 9:} For $d$ odd, if $A$ is a $d\times d$ symmetric matrix and $S=\tr(A)I-A$, where
$I$ is the $d\times d$ identity matrix, then
\begin{enumerate}[nolistsep]
\item[$1)$] If $A$ is of maximal trace over rotation matrices, then $S$ is positive
semidefinite.
\item[$2)$] $S$ fails to be positive semidefinite if and only if there exists a
rotation matrix $V$ such that $V=2vv^T-I$ for some vector $v\in \mathbb{R}^d$, $\|v\|=1$, and
$\tr(VA)>\tr(A)$.
\end{enumerate}
\,\smallskip
{\bf Proof:} Assume $A$ is of maximal trace over rotation matrices. If $S$ is not positive
semidefinite, then there is a vector $v\in \mathbb{R}^d$, $\|v\|=1$, such that
$$v^TSv=v^T(\tr(A)I-A)v<0.$$
Then $v^T\tr(A)Iv-v^TAv<0$ so that $v^TAv>\tr(A)v^Tv=\tr(A)$.\\
Let $V=2vv^T-I$. Then $-V$ is a Householder reflection matrix \cite{lay} which is well known to be
a symmetric orthogonal matrix of determinant equal to negative one.\\
Thus, as $d$ is odd it must be that det$(V)=1$ so that $V$ is a rotation matrix.\\
Note $VA=(2vv^T-I)A=2vv^TA-A$ so that
\begin{eqnarray*}
\tr(VA)&=&2\tr(vv^TA)-\tr(A)= 2\tr(v^TAv)-\tr(A)\\
&=&2v^TAv-\tr(A)>2\tr(A)-\tr(A)=\tr(A)
\end{eqnarray*}
which contradicts $A$ is of maximal trace over rotation matrices.
Thus, $S$ must be positive semidefinite and therefore statement~$1)$ holds.\\
From the proof of~$1)$ it is clear that~$2)$ holds.
\endproof
\,\medskip\par
We note the rest of this section is mostly concerned with results about $3\times 3$
matrices to be used in Section~5 for accomplishing the three-dimensional aspect of the secondary
goal of this paper: identifying alternative ways, other than the SVD,
of obtaining solutions to the problems of interest.\,\medskip\par
In what follows, when dealing with three-dimensional rotation matrices, given one such matrix,
say $W$, $W$ will be specified by an axis of rotation $w$, where $w$ is a unit vector in $\mathbb{R}^3$,
and a rotation angle $\theta$, $0\leq\theta\leq\pi$, where $\theta$ corresponds to a rotation
angle around the axis of rotation in a counterclockwise direction. The direction of the axis
of rotation $w$ is determined by the right-hand rule, i.e., the direction in which the thumb points
while curling the other fingers of the right hand around the axis of rotation with the curl
of the fingers representing a movement in the $\theta$ direction. Accordingly, given a $3\times 3$
rotation matrix $W$ with axis of rotation $w$ and rotation angle $\theta$ as just described,
assuming $w=(w_x,w_y,w_z)^T$, it is well known that
\[ W=\mathlarger{\left(\begin{smallmatrix}
\cos\theta+w_x^2(1-cos\theta) & w_xw_y(1-\cos\theta)-w_z\sin\theta
                                                            & w_xw_z(1-\cos\theta)+w_y\sin\theta\\
w_yw_x(1-cos\theta)+w_z\sin\theta & \cos\theta+w_y^2(1-\cos\theta)
                                                            & w_yw_z(1-\cos\theta)-w_x\sin\theta\\
w_zw_x(1-cos\theta)-w_y\sin\theta & w_zw_y(1-\cos\theta)+w_x\sin\theta
                                                            & \cos\theta+w_z^2(1-\cos\theta)\\
                     \end{smallmatrix}\right).}\]
Note $W=I$ for $\theta=0$, $I$ the $3\times 3$ identity matrix, $W=2ww^T-I$ for $\theta=\pi$, and
that given a $3\times 3$ symmetric matrix $A$
\[ A= \begin{pmatrix}
a_{11} & a_{12} & a_{13}\\
a_{21} & a_{22} & a_{23}\\
a_{31} & a_{32} & a_{33}\\
\end{pmatrix}\]
then it is not hard to show that
\begin{eqnarray*}
\tr(WA)&=&(a_{11}+a_{22}+a_{33})\cos\theta
+(a_{11}w_x^2+a_{22}w_y^2+a_{33}w_z^2\\
&&+2a_{12}w_xw_y+2a_{13}w_xw_z+2a_{23}w_yw_z)(1-\cos\theta)\\
&=& \tr(A)\cos\theta + w^TAw(1-\cos\theta).
\end{eqnarray*}
Thus, $\tr(WA)$ is an affine combination of $\tr(A)$ and $w^TAw$, where $\theta$ goes from
$0$ to~$\pi$. It follows then that $\tr(WA)$ achieves its minimum and maximum at either $\theta=0$
or $\theta=\pi$, and if it achieves its minimum (maximum) at $\theta=0$ then it must achieve
its maximum (minimum) at~$\theta=\pi$ and vice versa.
\,\medskip\par
Together with~$2)$ of Proposition 9, the following proposition provides another way of proving that
if $S$ of Proposition~8 is positive semidefinite, then matrix $A$ of the same proposition is of
maximal trace over rotation matrices.\,\medskip\\
{\bf Proposition 10:} If $A$ is a $3\times 3$ symmetric matrix and $W$ is any $3\times 3$ rotation
matrix with axis of rotation $w$ such that $\tr(WA)>\tr(A)$, then among all rotation matrices
$\hat{W}$ with axis of rotation~$w$, $\hat{W}=2ww^T-I$ maximizes $\tr(\hat{W}A)$ by a rotation
of $\pi$ radians. In particular, for this $\hat{W}$, $\tr(\hat{W}A)\geq\tr(WA)>\tr(A)$.
\,\medskip\\
{\bf Proof:} Because $\tr(WA)>\tr(A)$, then among all rotation matrices $\hat{W}$ with axis of
rotation $w$, $\hat{W}=I$ ($\theta=0$) must minimize $\tr(\hat{W}A)$ so that then
$\hat{W}=2ww^T-I$ ($\theta=\pi$) must maximize~$\tr(\hat{W}A)$.
\endproof
\,\medskip\par
Given a $3\times 3$ symmetric matrix $A$ that is not of maximal trace over rotation matrices,
the following proposition shows how to compute a $3\times 3$ rotation matrix $W$ such that
$WA$ is of maximal trace over rotation matrices if a unit eigenvector corresponding to the
largest eigenvalue of $A$ is~known.\,\medskip\\
{\bf Proposition 11:} Let $A$ be a $3\times 3$ symmetric matrix that is not of maximal trace
over rotation matrices. Let $\sigma$ be the largest eigenvalue of $A$, and $w$ a unit vector
in $\mathbb{R}^3$ that is an eigenvector of $A$ corresponding to $\sigma$. Let $W=2ww^T-I$.
Then $WA$ is of maximal trace over rotation matrices.
\,\medskip\\
{\bf Proof:} Let $V$ be any rotation matrix and let $v$ and $\theta$ be the rotation axis and
rotation angle associated with $V$, respectively.
Assume $\tr(VA)>\tr(A)$ and $\tr(VA)\geq\tr(\hat{V}A)$ for all rotation matrices
$\hat{V}$ with $v$ as the axis of rotation.
Then as above it must be that $\tr(VA) = \tr(A)\cos\theta + v^TAv(1-\cos\theta)$
with either $\theta=0$ or~$\theta=\pi$.\\
If $\theta=0$, then $\tr(VA)=\tr(A)$, a contradiction, thus it must be that~$\theta=\pi$
so that~$V=2vv^T-I$ and $\tr(VA)=-\tr(A)+2v^TAv$.\\
Accordingly, we look for a rotation matrix $W$ with axis of rotation $w$, such that
$\tr(WA)\geq\tr(VA)$ for all rotation matrices $V$, in particular any $V$ with
$\tr(VA)>\tr(A)$ and any $V$ with axis of rotation~$w$. Thus, if $W$ exists, it must be
that $W=2ww^T-I$, $v=w$ maximizing $v^TAv$.\\
Let $\sigma$ be the largest eigenvalue of $A$ and let $w$ be a unit eigenvector of $A$
corresponding to~$\sigma$. Then it is well known~\cite{lay} that $v=w$ maximizes~$v^TAv$
($\sigma$ the maximum value of $v^TAv$). Thus $W=2ww^T-I$ is as required.
\endproof
\,\medskip\par
Given a $d\times d$ symmetric matrix $A$, the following proposition shows how to compute a
$d\times d$ rotation matrix $W$ such that $WA$ is of maximal trace over rotation matrices
if an orthogonal diagonalization of $A$ is~known.\,\medskip\\
{\bf Proposition 12:} Let $A$ be a $d\times d$ symmetric matrix.
Let $V$, $D$ be $d\times d$ matrices such that $V$ is orthogonal,
$D=\mathrm{diag}\{\alpha_1,\ldots,\alpha_d\}$ with $\alpha_j$, \mbox{$j=1,\ldots,d$}, the eigenvalues
of $A$, and $A=V^TDV$. Define a set of integers~$H$ by
$$H=\{i\;|\;\alpha_i<0,i=1,\ldots,d\}.$$ If $H$ has an odd number of
elements, let $k= \arg\min_j\{|\alpha_j|,j=1,\ldots,d\}$. If $k\in H$ let $H=H\setminus\{k\}$.
Otherwise, let $H=H\cup\{k\}$. Let $G$ be the $d\times d$ orthogonal matrix with entries
$g_{lh}$, $l,h=1,\ldots,d$, among which, the nonzero entries are given by
$$g_{mm}=1,\ m=1,\ldots,d,\ m\not\in H,\ \ g_{mm}=-1,\ m=1,\ldots,d,\ m\in H.$$
Let $W=V^TGV$.  Then $W$ is a $d\times d$ rotation matrix and $WA$ is of maximal trace
over rotation matrices.\,\medskip\\
{\bf Proof:} Note det$(G)=1$ as $H$ is empty or has an even number of elements. Thus
det$(W)=1$ as well. Letting $\hat{D}=GD$, then $\hat{D}$ is a diagonal matrix with at most
one negative element in the diagonal, which, if it exists, is no larger in absolute value than
the other elements of the diagonal. Thus $V^T\hat{D}V$ must be of maximal trace over rotation
matrices. But $WA=V^TGVA=V^TGVV^TDV)=V^TGDV=V^T\hat{D}V$. Thus, $WA$ is of maximal trace
over rotation matrices.
\endproof
\section{\normalsize The Two-Dimensional Case: Computation without SVD}
In the two-dimensional case, it is possible to determine solutions to the problems of interest
in closed form that do not require the SVD method, i.e., the Kabsch-Umeyama algorithm. Suppose
$P=\{p_1,\ldots,p_n\}$, $Q=\{q_1,\ldots,q_n\}$
are each sets of $n$ points in~$\mathbb{R}^2$, and $W=\{w_1,\ldots,w_n\}$ is a set of $n$ nonnegative
numbers (weights). First we look at the problem of minimizing $\Delta(P,Q,U)$, i.e., of finding
a $2\times 2$ rotation matrix $U$ for which $\Delta(P,Q,U)$ is as small as possible. As we will
see, the problem of minimizing $\Delta(P,Q,W,U)$ can be approached in a similar manner with
some minor modifications. Here, for the sake of completeness, we first obtain the solutions
through a direct minimization of $\Delta(P,Q,U)$ and $\Delta(P,Q,W,U)$ that takes advantage
of various trigonometric identities
and of the representation of the points in terms of polar coordinates. However, as demonstrated
toward the end of this section, the trace maximization approach developed in Section~2 produces
the same solutions with a lot of less effort.
\par For each $i$, $i=1,\ldots,n$, let $p_i$ and $q_i$ be given in polar coordinates as
$p_i = (s_i, \sigma_i)$, $q_i = (r_i, \rho_i)$, where the first coordinate denotes the distance
from the point to the origin and the second denotes the angle (in radians) from the positive first
axis to the ray through the point from the origin. Clearly $s_i, r_i\geq 0$,
$0\leq \sigma_i,\rho_i<2\pi$, and if $p_i=(x_i, y_i)$, $q_i=(x_i^\prime, y_i^\prime)$,
in rectangular coordinates, then $x_i=s_i\cos\sigma_i$, $y_i=s_i\sin\sigma_i$,
$x_i^\prime=r_i\cos\rho_i$, $y_i^\prime=r_i\sin\rho_i$.
\,\smallskip\\
It is well known that if $U$ is a rotation matrix by $\theta$ radians in the counterclockwise
direction, $0\leq\theta<2\pi$, then
$$U=\mathlarger{\mathlarger{\left(\begin{smallmatrix}\cos\theta & \hspace{.05in} -\sin\theta\\
                     \sin\theta & \hspace{.15in} \cos\theta\\ \end{smallmatrix} \right).}}$$
Thus, using column vectors to perform the matrix multiplication
\begin{eqnarray*}
Uq_i&=&U\big(x_i^\prime,y_i^\prime\big)^T=\big(x_i^\prime\cos\theta-y_i^\prime\sin\theta,
                               x_i^\prime\sin\theta+y_i^\prime\cos\theta\big)^T\\
&=&\big(r_i\cos\rho_i\cos\theta-r_i\sin\rho_i\sin\theta,
                                      r_i\cos\rho_i\sin\theta+r_i\sin\rho_i\cos\theta\big)^T\\
&=&\big(r_i\cos(\rho_i+\theta),r_i\sin(\rho_i+\theta)\big)^T
\end{eqnarray*}
and
\begin{eqnarray*}
&& \hspace*{-.275in}\Delta(P,Q,U)=\sum_{i=1}^n||Uq_i - p_i||^2\\
&=& \sum_{i=1}^n\big((r_i \cos(\rho_i + \theta) - s_i\cos\sigma_i)^2
+(r_i \sin(\rho_i + \theta) - s_i\sin\sigma_i)^2\big) \\
&=& \sum_{i=1}^n\big(
r_i^2\cos^2(\rho_i+\theta)-2r_is_i\cos(\rho_i+\theta)\cos\sigma_i+s_i^2\cos^2\sigma_i\\
&&+\ r_i^2\sin^2(\rho_i+\theta)-2r_is_i\sin(\rho_i+\theta)\sin\sigma_i+s_i^2\sin^2\sigma_i\big)\\
&=&\sum_{i=1}^n \big(r_i^2 - 2r_is_i(\cos(\rho_i + \theta)\cos\sigma_i
+ \sin(\rho_i + \theta)\sin\sigma_i) + s_i^2\big)\\
&=&\sum_{i=1}^n \big(r_i^2 - 2r_is_i(\cos(\rho_i + \theta)\cos(-\sigma_i)
- \sin(\rho_i + \theta)\sin(-\sigma_i)) + s_i^2\big)\\
&=&\sum_{i=1}^n \big(r_i^2 - 2r_is_i\cos(\rho_i+\theta-\sigma_i) + s_i^2\big).
\end{eqnarray*}
Letting this last expression, which is equal to $\Delta(P,Q,U)$, be denoted by $f(\theta)$, then
\begin{eqnarray*}
f(\theta) &=& \sum_{i=1}^n \big(r_i^2 - 2r_is_i\cos(\rho_i-\sigma_i+\theta) + s_i^2\big)\\
&=&\sum_{i=1}^n \big(r_i^2 + s_i^2\big) - 2\sum_{i=1}^n
r_is_i\cos(\rho_i - \sigma_i + \theta) \\
&=&\sum_{i=1}^n \big(r_i^2 + s_i^2\big) - 2\sum_{i=1}^n
r_is_i\big(\cos(\rho_i - \sigma_i)\cos\theta-\sin(\rho_i-\sigma_i)\sin\theta\big)\\
&=& \sum_{i=1}^n\big(r_i^2 + s_i^2\big) - a\cos \theta +  b \sin \theta
\end{eqnarray*}
where
$$a =  2\sum_{i=1}^nr_is_i \cos(\rho_i-\sigma_i)$$
and
$$b = 2\sum_{i=1}^n r_i s_i \sin(\rho_i - \sigma_i).$$
Note that in terms of the rectangular coordinates of the points $p_i$, $q_i$
\begin{eqnarray*}
a &=& 2\sum_{i=1}^n r_i s_i \cos(\rho_i - \sigma_i)
            =2\sum_{i=1}^n r_i s_i(\cos\rho_i \cos \sigma_i
                 + \sin\rho_i \sin \sigma_i)  \\
            &=&2\sum_{i=1}^n (r_i \cos\rho_i s_i \cos \sigma_i
                 + r_i \sin\rho_i s_i \sin \sigma_i)
            =2\sum_{i=1}^n (x^\prime_i x_i + y^\prime_i y_i)\\
            &=&2\sum_{i=1}^n\,(x_i,y_i)\cdot(x^\prime_i,y^\prime_i)
\end{eqnarray*}
and
\begin{eqnarray*}
b &=& 2\sum_{i=1}^n r_i s_i \sin(\rho_i - \sigma_i)
            =2\sum_{i=1}^n r_i s_i(\sin \rho_i \cos \sigma_i
                 - \cos \rho_i \sin \sigma_i)  \\
            &=&2\sum_{i=1}^n (r_i \sin \rho_i s_i \cos \sigma_i
                 - r_i \cos \rho_i s_i \sin \sigma_i)
            =2\sum_{i=1}^n (y^\prime_i x_i - x^\prime_i y_i)\\
&=& 2\sum_{i=1}^n \mathlarger{\mathlarger{\left|\begin{smallmatrix}x_i & x^\prime_i\\
                     y_i & y^\prime_i\\ \end{smallmatrix} \right|.}}
\end{eqnarray*}
For each $i$, $i=1,\ldots,n$, letting $D_i$ be the dot product of $p_i$ and $q_i$,
then $a$ can be described as twice the sum of the $D_i$'s. On the other hand,
for each $i$, $i=1,\ldots,n$, letting $A_i$ be the signed area of the parallelogram
spanned by the vectors $\vec{0p_i}$, $\vec{0q_i}$, where the area is positive if
the angle in a counterclockwise direction from $\vec{0p_i}$ to $\vec{0q_i}$ is between
$0$ and $\pi$, zero or negative otherwise, then $b$ can be described as twice the
sum of the $A_i$'s.
\,\medskip\\
{\bf Theorem 1:}
Let $a=\sum_{i=1}^n x^\prime_i x_i + \sum_{i=1}^n y^\prime_i y_i$ and
$b=\sum_{i=1}^n y^\prime_i x_i - \sum_{i=1}^n x^\prime_i y_i$.
If $a = b = 0$, then $U=I$, $I$ the $2\times 2$ identity matrix,
minimizes~$\Delta(P,Q,U)$. Otherwise, with $c=\sqrt{a^2+b^2}$ and
\[\hat{U}= \mathlarger{\mathlarger{\left(\begin{smallmatrix} \hspace{.1in} a/c& \hspace{.05in} b/c\\
                     -b/c & \hspace{.05in} a/c\\ \end{smallmatrix} \right)}}\]
then $U=\hat{U}$ minimizes~$\Delta(P,Q,U)$.\,\medskip\\
{\bf Proof:} If $a = b = 0$, with $f$ as derived above, then
$f(\theta) = \sum_{i=1}^n \big(r_i^2 + s_i^2\big)$, i.e.,
$f(\theta)$ is constant so that $\Delta(P,Q,U)$ is constant as well, i.e., it has the same
value for all rotation matrices~$U$.
Thus, any $\theta$ minimizes $f(\theta)$, in particular $\theta=0$, and
therefore $U=I$, $I$ the $2\times 2$ identity matrix, minimizes~$\Delta(P,Q,U)$.\\
Otherwise, $f^\prime(\theta)=a\sin\theta+b\cos\theta$.\\
Since $a\,y+b\,x=0$ is the equation of a straight line $L$ through the origin, then $L$
must cross the unit circle at two points that are antipodal of each~other, and
it is easy to verify that these points are $(x,y)=(a/c,-b/c)$ and $(x,y)=(-a/c,b/c)$.\\
Since every point on the unit circle is of the form $(\cos\theta,\sin\theta)$ for some
$\theta$, $0\leq\theta <2\pi$, then for some $\theta_1$, $\theta_2$,
$0\leq\theta_1,\theta_2<2\pi$, it must be that $(a/c,-b/c)=(\cos\theta_1,\sin\theta_1)$
and $(-a/c,b/c)=(\cos\theta_2,\sin\theta_2)$.
Clearly, $f^\prime(\theta_1)=f^\prime(\theta_2)=0$.\\
Noting $f^{\prime\prime}(\theta)=a\cos\theta-b\sin\theta$, then
$f^{\prime\prime}(\theta_1)=a(a/c)-b(-b/c)=a^2/c+b^2/c>0$, and
$f^{\prime\prime}(\theta_2)=a(-a/c)-b(b/c)=-a^2/c-b^2/c<0$.\\
Thus, $f(\theta_1)$ is a local minimum of~$f$ on~$[0,2\pi)$ so that by the
differentiability and periodicity of~$f$ it is a global minimum of $f$ and,
therefore, $U=\hat{U}$ minimizes~$\Delta(P,Q,U)$.
\endproof
\,\medskip\par
With minor modifications due to the weights, arguing as above, a similar result can
be obtained for the more general problem. \,\medskip\\
{\bf Theorem 2:}
Let $a=\sum_{i=1}^n w_i x^\prime_i x_i + \sum_{i=1}^n w_i y^\prime_i y_i$ and
$b=\sum_{i=1}^n w_i y^\prime_i x_i - \sum_{i=1}^n w_i x^\prime_i y_i$.
If $a = b = 0$, then $U=I$, $I$ the $2\times 2$ identity matrix,
minimizes~$\Delta(P,Q,W,U)$. Otherwise, with $c=\sqrt{a^2+b^2}$ and
\[\hat{U}=\mathlarger{\mathlarger{\left(\begin{smallmatrix} \hspace{.1in} a/c& \hspace{.05in} b/c\\
                     -b/c & \hspace{.05in} a/c\\ \end{smallmatrix} \right)}}\]
then $U=\hat{U}$ minimizes~$\Delta(P,Q,W,U)$.\,\medskip\\
{\bf Proof:} The same as that of Theorem~1 with minor modifications.
\endproof
\, \medskip\par
Finally, let
\[a_{11}=\sum_{i=1}^n w_i x^\prime_i x_i,\ a_{22}=\sum_{i=1}^n w_i y^\prime_i y_i,
\ a_{21}=\sum_{i=1}^n w_i y^\prime_i x_i,\ a_{12}=\sum_{i=1}^n w_i x^\prime_i y_i\]
and note with $a$ and $b$ as above that $a=a_{11}+a_{22}$, $b=a_{21}-a_{12}$.\\
If
\[A=\mathlarger{\mathlarger{\left(\begin{smallmatrix} a_{11} & \hspace{.05in} a_{12}\\
                     a_{21} & \hspace{.05in} a_{22}\\ \end{smallmatrix} \right)}}\]
then minimizing $\Delta(P,Q,W,U)$ is equivalent, as observed in Section~2, to maximizing
$\tr(UA)$ over all $2\times 2$ rotation matrices~$U$, where if $U$ is a rotation matrix by
$\theta$ radians in a counterclockwise direction, $0\leq\theta<2\pi$, then
\[ U=\mathlarger{\mathlarger{\left(\begin{smallmatrix}\cos\theta & \hspace{.05in} -\sin\theta\\
                     \sin\theta & \hspace{.15in} \cos\theta\\ \end{smallmatrix} \right).}}\]
Note then that\\
\hspace*{.15in} $\tr(UA) = (a_{11}\cos\theta-a_{21}\sin\theta)+(a_{12}\sin\theta+a_{22}\cos\theta)$\\
\hspace*{.715in} $=(a_{11}+a_{22})\cos\theta -(a_{21}-a_{12})\sin\theta=a\cos\theta-b\sin\theta$\\
so that by using the trace maximization approach, we have essentially derived the
function $f$, previously derived above, with a lot of less effort.\\
Note that if $a=b=0$, then clearly $a_{11}=-a_{22}$ and~$a_{21}=a_{12}$.
Also as established above it must be that $\Delta(P,Q,W,U)$ has the same
value for all rotation matrices $U$, and, therefore, so does~$\tr(UA)$.\\
Thus, it is no coincidence that given any arbitrary $\theta$, $0\leq \theta<2\pi$,
if $U$ is the rotation matrix by $\theta$ radians in a counterclockwise direction, then\\
$\tr(UA) = (a_{11}\cos\theta-a_{21}\sin\theta)+(a_{21}\sin\theta-a_{11}\cos\theta) = 0$, i.e.,
$\tr(UA)=0$ for all rotation matrices~$U$.\\
Also $(a_{11}\sin\theta+a_{21}\cos\theta)-(a_{21}\cos\theta+a_{11}\sin\theta) = 0$,
so that $UA$ is indeed a symmetric matrix.\\
On the other hand, if $a\not=0$ or $b\not=0$, with $c=\sqrt{a^2+b^2}$, then $U=\hat{U}$ that
minimizes $\Delta(P,Q,W,U)$ in Theorem~2 above must maximize~$\tr(UA)$ and the maximum is\\
$\tr(\hat{U}A)=a_{11}(a/c)-a_{21}(-b/c)+a_{12}(-b/c)+a_{22}(a/c)$\\
$=(a_{11}+a_{22})(a/c)+(a_{21}-a_{12})(b/c)=a(a/c)+b(b/c)= (a^2+b^2)/c=c>0$ which is
nonnegative, actually positive, as expected according to Corollary 1 of Section~3.\\
We also have the relation $a_{11}(-b/c)+a_{21}(a/c)-a_{12}(a/c)+a_{22}(-b/c)$\\
$=(a_{11}+a_{22})(-b/c)+(a_{21}-a_{12})(a/c)=-ab/c+ba/c=0$, so that $\hat{U}A$ is indeed
a symmetric matrix.
\section{\normalsize The Three-Dimensional Case: Computation without SVD}
Given a real $3\times 3$ matrix $M$ that is not of maximal trace over rotation matrices,
in this section, if the matrix $M$ is symmetric, we present an approach that does not use
the SVD method, i.e., the Kabsch-Umeyama algorithm, for computing a $3\times 3$ rotation
matrix $U$ such that $UM$ is of maximal trace over rotation matrices. This approach, which
is based on a trigonometric identity, is a consequence of Proposition~$11$ in Section~3, and
if the matrix $M$ is not symmetric, part of it can still be used to produce the usual
orthogonal matrices necessary to carry out the SVD method. Being able to find such a
matrix $U$ for a matrix $M$, not necessarily symmetric, is what is required to solve Wahba's
problem, not only for $3\times 3$ matrices, but also for $d\times d$ matrices
for any~$d$,~$d\geq 2$. As described in Section~1 and Section~2 of this paper, in Wahba's problem
the number $\Delta(P,Q,W,U)$ is minimized, where $P=\{p_i\}$, $Q=\{q_i\}$, $i=1,\ldots,n$,
are each sets of $n$ points in $\mathbb{R}^d$, and $W=\{w_i\}$, $i=1,\ldots,n$, is a set of $n$
nonnegative weights. Accordingly, in the three-dimensional version of the problem,
the points $p_i$, $q_i$ are then of the form $p_i=(x_i,y_i,z_i)$, $q_i=(x_i',y_i',z_i')$,
$i=1,\ldots,n$, and with
\begin{eqnarray*}
  m_{11}=\sum_{i=1}^n w_i x^\prime_i x_i, \ m_{12}=\sum_{i=1}^n w_i x^\prime_i y_i,
\ m_{13}=\sum_{i=1}^n w_i x^\prime_i z_i\\
  m_{21}=\sum_{i=1}^n w_i y^\prime_i x_i, \ m_{22}=\sum_{i=1}^n w_i y^\prime_i y_i,
\ m_{23}=\sum_{i=1}^n w_i y^\prime_i z_i\\
  m_{31}=\sum_{i=1}^n w_i z^\prime_i x_i, \ m_{32}=\sum_{i=1}^n w_i z^\prime_i y_i,
\ m_{33}=\sum_{i=1}^n w_i z^\prime_i z_i
\end{eqnarray*}
the $3\times 3$ matrix of interest $M$ is then
\[ M= \begin{pmatrix}
m_{11} & m_{12} & m_{13}\\
m_{21} & m_{22} & m_{23}\\
m_{31} & m_{32} & m_{33}\\
\end{pmatrix}.\]
If $M$ is a symmetric matrix, in this approach we refer to $M$ by the name $A$ to signify that
$A\,(=M)$ is symmetric, and if $A$ is not of maximal trace over rotation matrices, a $3\times 3$
rotation matrix $R$ is computed without using the SVD method in such a way that $RA$ is of maximal
trace over rotation matrices. We note that before trying to compute $R$, the matrix $A$ should be
tested for the maximality of the trace. This can be done as a consequence of Proposition~$8$ in
Section~3, i.e., by testing whether $S=\tr(A)I-A$ is positive semidefinite, where $I$ is the
$3\times 3$ identity matrix. It is well known that a square matrix is positive semidefinite if and
only if all its principal minors are nonnegative. Since positive definiteness implies positive
semidefiniteness, and because a square matrix is positive definite if and only if
all its leading principal minors are positive, we test the matrix first for
positive definiteness as a $3\times 3$ matrix has seven principal minors of which only
three are of the leading kind. On the other hand, if the matrix $M$ is not symmetric,
part of the approach can still be used on $A=M^TM$ which is symmetric, to produce the usual
orthogonal matrices necessary to carry out the SVD method.
\,\medskip\par
The approach which we present next is a consequence of Proposition~$11$ in Section~3. According
to the proposition, if $A\,(=M)$ is a $3\times 3$ symmetric matrix that is not of maximal trace over
rotation matrices, then in order to obtain a $3\times 3$ rotation matrix $R$ so that $RA$ is
of maximal trace over rotation matrices, it suffices to compute $R=2\hat{r}\,\hat{r}^T-I$,
where $\hat{r}$ is a unit vector in $\mathbb{R}^3$ that is an eigenvector of~$A$ corresponding
to the largest eigenvalue of~$A$. In our approach, the computation of $\hat{r}$, and, if
necessary, the computation of all eigenvectors of $A=M^TM$ (to carry out the SVD method if
$M$ is not symmetric), is essentially as presented in~\cite{eberly,smith}. We note that a nice
alternative method can be found in~\cite{mladenov} which is a two-step procedure based on a
vector parametrization of the group of three-dimensional rotations.
Following ideas in~\cite{eberly,smith}, we accomplish our purpose
by taking advantage of a $3\times 3$ matrix $B$ that is a linear combination of $A$
and $I$ in the appropriate manner so that the characteristic polynomial of $B$ is such
that it allows the application of a trigonometric identity in order to obtain its roots
in closed form and thus those of the characteristic polynomial of~$A$. \smallskip\\
Thus, let $A$ be a $3\times 3$ symmetric matrix (we do not assume $A$ is not of
maximal trace over rotation matrices at this point)
\[ A= \begin{pmatrix}
a_{11} & a_{12} & a_{13}\\
a_{21} & a_{22} & a_{23}\\
a_{31} & a_{32} & a_{33}\\
\end{pmatrix}.\]
It is well known that if $A$ is just any $3\times 3$ matrix, the characteristic polynomial
of $A$ is
\[f(\alpha) = \mathrm{det}(\alpha I-A)
= \alpha^3-\alpha^2\tr(A)-\alpha 1/2(\tr(A^2)-\tr^2(A))-\mathrm{det}(A).\]
Given numbers $p>0$ and $q$, we define a $3\times 3$ matrix $B$ by $B=(A-qI)/p$ so that
$A=pB+Iq$. Note that if $v$ is an eigenvector of $A$ corresponding to an eigenvalue $\alpha$
of $A$, i.e., $Av=\alpha v$, then $Bv=((\alpha-q)/p)v$ so that $v$ is an eigenvector of $B$
corresponding to the eigenvalue $(\alpha-q)/p$ of~$B$. Conversely, if $v$ is an eigenvector of
$B$ corresponding to an eigenvalue $\beta$ of $B$, i.e., $Bv=\beta v$, then $Av = (p\beta +q)v$
so that $v$ is an eigenvector of $A$ corresponding to the eigenvalue $p\beta +q$ of~$A$.
Thus, $A$ and $B$ have the same~eigenvectors.\\
Let $q=\tr(A)/3$ and $p=(\tr((A-qI)^2)/6)^{1/2}$. Then $p\geq0$.\\
We treat $p=0$ as a special case so that then we can assume $p>0$ as required.\\
Accordingly, we note that $p=0$ if and only if $\tr((A-qI)^2) = 0$, and since it is readily shown
that $A$ is a symmetric matrix, then
\[(a_{11}-q)^2+(a_{22}-q)^2+(a_{33}-q)^2+2a_{12}^2+2a_{13}^2+2a_{23}^2=0.\]
Thus $a_{11}=a_{22}=a_{33}=q$ and $a_{12}=a_{21}=a_{13}=a_{31}=a_{23}=a_{32}=0$ so that
$A=\mathrm{diag}\{q,q,q\}$ and $q$ is therefore the only eigenvalue (a multiple eigenvalue) of~$A$.\\
Assuming now $p>0$ as required, then, in particular, $\tr((A-qI)^2)\not=0$.\\ Note\\
$\tr(B)=\tr((A-qI)/p)= 1/p(\tr(A)-\tr(qI))=1/p(\tr(A)-3(\tr(A)/3))=0$\\ and\\
$\tr(B^2)=\tr(((A-qI)/p)^2)=\tr(((A-qI)/((\tr((A-qI)^2)/6)^{1/2}))^2)$\\
$\hspace*{.475in}=\tr((A-qI)^2/(\tr((A-qI)^2)/6))$\\
$\hspace*{.475in}=(6/(\tr((A-qI)^2)))(\tr((A-qI)^2))=6$.
\smallskip\\
Thus, the characteristic polynomial of $B$ is
\[g(\beta) = \mathrm{det}(\beta I-B)=\beta^3-3\beta-\mathrm{det}(B).\]
We show $|\mathrm{det}(B)|\leq 2$. The general cubic equation has the form $ax^3+bx^2+cx+d=0$
with $a\not= 0$. It is well known that the numbers of real and complex roots are determined by
the discriminant $\Delta$ of this equation, $\Delta =18abcd-4b^3d+b^2c^2-4ac^3-27a^2d^2$.\\
\hspace*{.1in} If $\Delta > 0$, then the equation has three distinct real roots.\\
\hspace*{.1in} If $\Delta = 0$, then it has a multiple root and all of its roots are real.\\
\hspace*{.1in} If $\Delta < 0$, then it has one real root and two complex conjugate roots.\\
For $g$ above, $a=1$, $b=0$, $c=-3$, $d=-\mathrm{det}(B)$, and since $B$ is clearly a symmetric
matrix, then it has three real roots.\\
Thus $\Delta=-4(-3)^3-27(-\mathrm{det}(B))^2=4\cdot27-27(\mathrm{det}(B))^2\geq 0$\\
so that $(\mathrm{det}(B))^2\leq 4$ and $|\mathrm{det}(B)|\leq 2$.\,\smallskip\\
Note the first derivative of $g$ is $g^\prime(\beta)=3\beta^2-3$ and $g^\prime(\beta)=0$ at
$\beta=-1$ and $\beta=1$. The second derivative is $g^{\prime\prime}(\beta)=3\beta$ and
$g^{\prime\prime}(-1)=-3$, $g^{\prime\prime}(1)=3$, so that $g$ has a local maximum at $\beta=-1$
and a local minimum at $\beta=1$.\\
Note as well $g(-2)=g(1)=-2-\mathrm{det}(B)$, $g(-1)=g(2)=2-\mathrm{det}(B)$ so that it is
not hard to see that for $-2<\mathrm{det}(B)<2$, $g$ alternates enough between positive
and negative values to have three distinct roots as predicted by its discriminant, all
in the interval~$(-2,2)$. Similarly, for $\mathrm{det}(B)=-2$ and $\mathrm{det}(B)=2$, it is not
hard to see that $g$ has two roots, one multiple, also as predicted by its discriminant, both in
the interval~$[-2,2]$.\\
Let $\beta_1$, $\beta_2$, $\beta_3$ be the three roots of $g$,
$-2\leq\beta_1\leq\beta_2\leq\beta_3\leq2$.\\
For $\theta\in [0,\pi]$, define $h:[0,\pi]\rightarrow [-2,2]$ by $h(\theta)=2\cos\theta$.
Clearly $h$ is one-to-one and onto, $h(0)=2$, $h(\pi)=-2$, so that numbers $\theta_1$,
$\theta_2$, $\theta_3$ exist such that $\pi\geq\theta_1\geq\theta_2\geq\theta_3\geq 0$,
$h(\theta_1)=\beta_1$, $h(\theta_2)=\beta_2$, $h(\theta_3)=\beta_3$.\\
Thus, $g(h(\theta))=(2\cos\theta)^3-3(2\cos\theta)-\mathrm{det}(B)=
2(4\cos^3\theta-3\cos\theta)-\mathrm{det}(B)$ has roots $\theta_1$, $\theta_2$, $\theta_3$
as just described, and since $\cos 3\theta=4\cos^3\theta-3\cos\theta$, then
$g(h(\theta))=2\cos 3\theta-\mathrm{det}(B)$ so that $\cos 3\theta=\mathrm{det}(B)/2$
at $\theta=\theta_1,\theta_2,\theta_3$. As then it must be that
$0\leq 3\theta_3\leq\pi\leq 3\theta_2\leq 2\pi\leq 3\theta_1\leq 3\pi$, it follows that\\
$3\theta_3=\arccos(\mathrm{det}(B)/2)$, $3\theta_2=2\pi-3\theta_3$, $3\theta_1=2\pi+3\theta_3$.\\
Thus, $\theta_3=\arccos(\mathrm{det}(B)/2)/3$, $\theta_2=2\pi/3-\theta_3$,
$\theta_1=2\pi/3+\theta_3$, from which $\beta_k$, $\alpha_k$, $k=1,2,3$, the eigenvalues of $B$
and $A$, respectively, can be computed as $\beta_k=2\cos\theta_k$,
and $\alpha_k=p\beta_k+q$, $k=1,2,3$. As $\beta_3$ is the largest eigenvalue of $B$ and
since $p>0$, then $\alpha_3$ must be the largest eigenvalue of~$A$.
\,\smallskip\\ From this discussion the next theorem follows.
\,\medskip\\
{\bf Theorem 3:} Let $A$ be a $3\times 3$ symmetric matrix. Let $a_{ij},i,j=1,2,3$, be the entries
of~$A$. With $I$ the $3\times 3$ identity matrix, let $q=\tr(A)/3=(a_{11}+a_{22}+a_{33})/3$,
$p=(\tr((A-qI)^2)/6)^{1/2}=
(((a_{11}-q)^2+(a_{22}-q)^2+(a_{33}-q)^2+2a_{12}^2+2a_{13}^2+2a_{23}^2)/6)^{1/2}$. If $p=0$, then
$A=\mathrm{diag}\{q,q,q\}$ and letting $\alpha_k=q$, $k=1,2,3$, then $\alpha_k$, $k=1,2,3$, are
the eigenvalues (the same eigenvalue) of~$A$.
Otherwise, let $B=(A-qI)/p$. Let $\theta_3=\arccos(\mathrm{det}(B)/2)/3$,
$\theta_2=2\pi/3-\theta_3$, $\theta_1=2\pi/3+\theta_3$. Let $\alpha_k=2p\cos\theta_k+q$, $k=1,2,3$.
Then $\alpha_1\leq\alpha_2\leq\alpha_3$, and $\alpha_k$, $k=1,2,3$, are the eigenvalues of~$A$.
\,\medskip\\
Finally, given an eigenvalue $\alpha$ of $A$, $A$ a real $3\times 3$ symmetric matrix,
we show how to compute an orthonormal set of eigenvectors of $A$ that spans the eigenspace
of $A$ corresponding to~$\alpha$.\\
For this purpose, let $C=A-\alpha I$. If $C$ is the zero matrix, then
$A=\mathrm{diag}\{\alpha,\alpha,\alpha\}$ which incidentally, if $A$ is not of maximal
trace over rotation matrices, can only happen if $\alpha<0$ by Corollary~1 of Section~3.
Since any vector in $\mathbb{R}^3$ is then an eigenvector of $A$ corresponding to the only
eigenvalue $\alpha$ of $A$, then, for example, $\{(1,0,0)^T,(0,1,0)^T,(0,0,1)^T\}$ is
an orthonormal set of eigenvectors of $A$ that spans the eigenspace of $A$ corresponding
to~$\alpha$ (the eigenspace is all of~$\mathbb{R}^3$).\\
Thus, we assume $C$ is not the zero matrix so that the null space of $C$ is not all of
$\mathbb{R}^3$, and we already know, since $\alpha$ is an eigenvalue of $A$, that the null space of~$C$
does not consist exactly of the single point~$(0,0,0)^T$. Thus, the dimension of the null space
of $C$ is either one or two. As $C$ is clearly a symmetric matrix then its null space is the
orthogonal complement of its column space and the dimension of its column space, therefore,
can only be one or two as well.\\
Let $c_1$, $c_2$, $c_3$ be the column vectors of~$C$, and with $\times$ denoting
the cross product operation, let $v_1=c_1\times c_2$, $v_2=c_2\times c_3$, $v_3=c_3\times c_1$.
If one or more of the vectors $v_1$, $v_2$, $v_3$, is not zero, i.e., is not~$(0,0,0)^T$,
let $v$ be one such vector. Then $\|v\|\not=0$ and the two column vectors of $C$ whose cross
product is~$v$ span the column space of~$C$ (the dimension of the column space of~$C$ equals two
so that the dimension of the null space of~$C$ is one). Since $v$ is orthogonal to both, it must be
that $v$ is in the null space of~$C$ and $\hat{v}=v/\|v\|$ is then a unit vector that is an
eigenvector of~$A$ corresponding to the eigenvalue~$\alpha$ of~$A$. Thus, $\{\hat{v}\}$ is an
orthonormal set of eigenvectors of $A$ (one eigenvector) that spans the eigenspace of $A$
corresponding to~$\alpha$.\\
Finally, if all of $v_1$, $v_2$, $v_3$ equal~$(0,0,0)^T$,
then the dimension of the column space of~$C$ equals one so that the dimension of the null space
of~$C$ is two, and one or more of the column vectors
$c_1$, $c_2$, $c_3$, is not zero, i.e., is not~$(0,0,0)^T$. Let $u$ be one such vector and
let $w=(1,1,1)^T$. Clearly $u$ spans the column space of~$C$. With $u=(u_1,u_2,u_3)^T$, let
$k=\arg\,\max_j\{|u_j|, j=1,2,3\}$.  In the vector~$w$ replace the $k^{th}$ coordinate with~$0$.
Then $v_1=u\times w$ is not~$(0,0,0)^T$. Thus $\|v_1\|\not=0$, $v_1$ is orthogonal to $u$, and it
must be that $v_1$ is in the null space of~$C$ and $\hat{v}_1=v_1/\|v_1\|$ is then a unit vector that
is an eigenvector of~$A$ corresponding to the eigenvalue~$\alpha$ of~$A$. Furthermore,
$v_2=v_1\times u$ is not~$(0,0,0)^T$, thus $\|v_2\|\not=0$, $v_2$ is orthogonal to $v_1$ and to
$u$, and it must be that $v_2$ is also in the null space of~$C$. It follows then that
$\hat{v}_2=v_2/\|v_2\|$ is a unit vector that is an eigenvector of~$A$ orthogonal to~$\hat{v}_1$
corresponding to the eigenvalue~$\alpha$ of~$A$. Thus, $\{\hat{v}_1,\hat{v}_2\}$ is an orthonormal
set of eigenvectors of $A$ that spans the eigenspace of $A$ corresponding to~$\alpha$ ($\alpha$ is
of multiplicity two). Note that in this case we can actually identify a third eigenvector $\hat{u}$
of $A$ of unit length corresponding to the eigenvalue of~$A$ not equal to $\alpha$ by setting
$\hat{u}=u/\|u\|$.\,\smallskip\\
From this discussion the next theorem follows.\,\medskip\\
{\bf Theorem 4:} Let $A$ be a $3\times 3$ symmetric matrix. Let $\alpha$ be an eigenvalue of~$A$.
Let $C=A-\alpha I$. If $C$ is the zero matrix, then $\{(1,0,0)^T,(0,1,0)^T,(0,0,1)^T\}$
is an orthonormal set of eigenvectors of $A$ that spans the eigenspace of $A$ corresponding
to~$\alpha$ (the eigenspace is all of~$\mathbb{R}^3$). Otherwise, let $c_1$, $c_2$, $c_3$ be the column
vectors of~$C$, and let $v_1=c_1\times c_2$, $v_2=c_2\times c_3$, $v_3=c_3\times c_1$.
If one or more of the vectors $v_1$, $v_2$, $v_3$, is not zero, i.e., is not~$(0,0,0)^T$,
let $v$ be one such vector. Let $\hat{v}=v/\|v\|$. Then $\{\hat{v}\}$ is an orthonormal
set of eigenvectors of $A$ (one eigenvector) that spans the eigenspace of $A$ corresponding
to~$\alpha$. Otherwise, if all of $v_1$, $v_2$, $v_3$ equal~$(0,0,0)^T$, let $u$
be one of $c_1$, $c_2$, $c_3$, that is not zero, i.e., is not~$(0,0,0)^T$, and let
$w=(1,1,1)^T$.  With $u=(u_1,u_2,u_3)^T$, let $k=\arg\,\max_j\{|u_j|, j=1,2,3\}$.
In the vector~$w$ replace the $k^{th}$ coordinate with~$0$ and let $v_1=u\times w$.
Let $\hat{v}_1=v_1/\|v_1\|$. Furthermore, let $v_2=v_1\times u$ and $\hat{v}_2=v_2/\|v_2\|$.
Then $\{\hat{v}_1,\hat{v}_2\}$ is an orthonormal set of eigenvectors of $A$ that spans the eigenspace
of $A$ corresponding to~$\alpha$ ($\alpha$ is of multiplicity two).
\,\medskip\\
Given a real $3\times 3$ symmetric matrix $A$ that is not of maximal trace over rotation
matrices, then $\alpha_3$, of computation as described in Theorem~3, is the largest eigenvalue
of~$A$. Let $\hat{r}$ be any unit eigenvector of $A$ corresponding to the eigenvalue $\alpha_3$
of~$A$, of computation as described in Theorem~4 with $\alpha=\alpha_3$. Then, by
Proposition~$11$ in Section~3, if $R=2\hat{r}\,\hat{r}^T-I$, then $RA$ is of maximal trace over
rotation~matrices. On the other hand, given a real $3\times 3$ matrix $M$ that is not symmetric,
letting $A=M^TM$, then $A$ is symmetric and it is $A$ that is usually used to compute the SVD of~$M$.
Accordingly, Theorem~3 and Theorem~4 can then be used to compute an orthonormal basis of $\mathbb{R}^3$
consisting of eigenvectors of $A$ in the proper order which are then used to produce the
usual orthogonal matrices necessary to carry out the SVD method~\cite{lay}.
We note that all of the above has been successfully implemented in~Fortran.
Links to the code are provided in the next section.
\section{\normalsize The Three-Dimensional Case Revisited}
In this section, given a $3\times 3$ matrix $M$ that is not symmetric we describe a procedure
that uses the so-called Cayley transform \cite{bellman,cayley,weyl} in conjunction with Newton's
method to find a $3\times 3$ rotation matrix $U$ so that $UM$ is symmetric, possibly of maximal
trace over rotation matrices. If the resulting $UM$ is not of maximal trace over rotation matrices,
using the fact that $UM$ is symmetric, another $3\times 3$ rotation matrix $R$ can then be computed
(without the SVD) as described in the previous section so that $RUM$ is of maximal trace over
rotation matrices. Since the possibility exists that Newton's method can fail, whenever this
occurs the SVD method is carried out as just described at the end of the previous section.
\,\medskip\par Given a $d\times d$ matrix $B$ such that $I+B$ is invertible, $I$ the $d\times d$ identity
matrix, we denote by $C(B)$ the $d\times d$ matrix $$C(B) = (I-B)(I+B)^{-1}.$$ The matrix $C(B)$ is
called the {\em Cayley transform} of $B$ and it is well known \cite{bellman,cayley,weyl} that if
$C(B)$ exists, then $I+C(B)$ is invertible so that $C(C(B))$ exists and it is actually equal to~$B$.
\,\medskip\par Letting $A$ be any $d\times d$ skew-symmetric matrix ($A^T=-A$), then it is well known
\cite{bellman,cayley,weyl} that $I+A$ is invertible, and $Q=C(A)$ is a rotation matrix ($Q^TQ=I$,
det$(Q)=1)$). Conversely, letting $Q$ be any $d\times d$ orthogonal matrix with $I+Q$ invertible,
i.e., $-1$ is not an eigenvalue of $Q$, then it is also well known that $A=C(Q)$ is skew-symmetric.
Note that $-1$ not being an eigenvalue of $Q$ excludes at least all orthogonal matrices of
determinant negative one. In particular, for $d=3$, among rotation matrices, it excludes exactly all
rotation matrices whose rotation angle equals $\pi$ radians. Consequently, from the
above comments, for every $d\times d$ rotation matrix $Q$ with $I+Q$ invertible, there is a
$d\times d$ skew-symmetric matrix $A$ with $C(A)=Q$, and for every $d\times d$ rotation matrix
$Q$ with $I+Q$ not invertible there is no $d\times d$ skew-symmetric matrix $A$ with~$C(A)=Q$.
\,\medskip\\Given a $3\times 3$ skew-symmetric matrix $A$
$$A = \begin{pmatrix}\hspace{.1in} 0 & \hspace{.1in} r & -s \hspace{.03in} \cr
               -r & \hspace{.1in} 0 & \hspace{.1in} t \cr
               \hspace{.1in} s & -t & \hspace{.1in} 0 \cr
\end{pmatrix}
$$
then with $\Delta = 1 + r^2 + s^2 + t^2$ it is well known that
\begin{eqnarray*}
&& \frac{\Delta}{2}C(A)=\frac{\Delta}{2}I-A+A^2=\\
&& \begin{pmatrix}\frac{\Delta}{2} & 0 & 0 \cr
               0& \frac{\Delta}{2} & 0 \cr
               0 & 0& \frac{\Delta}{2}\cr\end{pmatrix} -
               \begin{pmatrix}\hspace{.1in} 0 & \hspace{.1in} r & -s \hspace{.03in} \cr
               -r & \hspace{.1in} 0 & \hspace{.1in} t \cr
               \hspace{.1in} s & -t & \hspace{.1in} 0 \cr \end{pmatrix}
+\begin{pmatrix}-r^2-s^2 & st & rt \cr
               st & -r^2 - t^2 & rs \cr rt & rs & -s^2 - t^2 \cr\end{pmatrix}.
\end{eqnarray*}
\par As we have seen, given a $d\times d$ matrix $M$ and a $d\times d$ rotation matrix~$U$, a
necessary condition for $UM$ to be of maximal trace over rotation matrices is that $UM$ be symmetric.
For~$d=3$ we use Newton's method in the way described below on some function $g$, defined below, in
order to find a rotation matrix $U$ such that $UM$ is symmetric, by finding a zero of~$g$. Since $U$
exists for which $UM$ is of maximal trace over rotation matrices, we know such a $U$ exists.
However, since the way in which we use Newton's method below is based on the Cayley transform,
if $U$ for which $UM$ is of maximal trace over rotation matrices is a rotation matrix whose
rotation angle equals~$\pi$, then Newton's method could
run into difficulties (the function $g$ on which Newton's method is used may not have a zero).
That $U$, $M$ do exist where $U$ is a rotation matrix whose rotation angle equals~$\pi$, $M$
is not symmetric, and $UM$ is of maximal trace over rotation matrices, is exemplified as follows
$$U = \begin{pmatrix}-1 & \hspace{.1in} 0 & \hspace{.1in} 0 \hspace{.03in} \cr
               \hspace{.1in} 0 & -1 & \hspace{.1in} 0 \cr
               \hspace{.1in} 0 & \hspace{.1in} 0 & \hspace{.1in} 1 \cr\end{pmatrix},
\ \ \ \ \ M = \begin{pmatrix}-2 & -1 & \hspace{.1in} 0 \hspace{.03in} \cr
               -1 & -2 & -1 \cr
                \hspace{.1in} 0 & \hspace{.1in} 1 & \hspace{.1in} 2 \cr\end{pmatrix},
\ \ \ \ \ UM = \begin{pmatrix}2 & 1 & 0 \cr
               1 & 2 & 1 \cr
                0 & 1 & 2 \cr\end{pmatrix}.
$$
\par Given $x = (r, s, t)^T \in \mathbb{R}^3$, let $A(x)$ denote
the matrix $A$ above.  Let
$$F(x) = \frac{(1 + r^2 + s^2 + t^2)}{2}C(A(x)).$$
Finally, define $g:\mathbb{R}^3\rightarrow \mathbb{R}^3$ by
$$g(x) =  (u, v, w)^T$$
where
$$
F(x)M-M^T F(x)^T = \begin{pmatrix} \hspace{.1in} 0& \hspace{.1in} u& -v \hspace{.03in} \cr
-u& \hspace{.1in} 0 & \hspace{.1in} w\cr \hspace{.1in} v&-w& \hspace{.1in} 0\end{pmatrix}.
$$
We wish to find a zero $\bar x = (\bar r, \bar s, \bar t)^T $ of $g$, i.e, $\bar{x}$ such that
$g(\bar x) = (0,0,0)$. Clearly, if $\delta= 1 + \bar{r}^2 + \bar{s}^2 + \bar{t}^2$, then
$U=\frac{2}{\delta}F(\bar{x})$ is a rotation matrix such that $UM$ is symmetric.
For this purpose we use Newton's method on~$g$.
\par Newton's method consists of performing a sequence of iterations based on the function $g$
and its Jacobian matrix~$J$, beginning from an initial point $x_0 \in \mathbb{R}^3$
\begin{eqnarray*}
x_0 &=& \mathrm{initial\ point\ in}\ \mathbb{R}^3\\
x_{k+1} &=& x_k -J(x_k)^{-1}g(x_k)\ \mathrm{for}\ k=0,1,2,\ldots
\end{eqnarray*}
Given that $g$ is sufficiently smooth and the Jacobian $J$ of $g$ is nonsingular at each~$x_k$,
if the initial point $x_0$ is ``sufficiently'' close to a root $\bar x$ of~$g$, then the sequence
$\{x_k\}$ converges to $\bar x$ and the rate of convergence is quadratic. Clearly, besides the
situation mentioned above, Newton's method could also run into difficulties if the initial point
$x_0$ is not close enough to a root of~$g$ or if the Jacobian of~$g$ is singular at some~$x_k$.
\par With $F$, $x$ and $\Delta$ as above, then again
\begin{eqnarray*}
&&F(x) = \frac{1 + r^2 + s^2 + t^2}{2} I - A + A^2=\\
&& \begin{pmatrix}\frac{\Delta}{2} & 0 & 0 \cr
               0& \frac{\Delta}{2} & 0 \cr
               0 & 0& \frac{\Delta}{2}\cr\end{pmatrix} -
               \begin{pmatrix} \hspace{.1in} 0 & \hspace{.1in} r & -s \hspace{.03in} \cr
               -r & \hspace{.1in} 0 & \hspace{.1in} t  \cr
               \hspace{.1in} s & -t & \hspace{.1in} 0 \cr\end{pmatrix}
 + \begin{pmatrix}-r^2-s^2 & st & rt \cr
               st & -r^2 - t^2 & rs \cr rt & rs & -s^2 - t^2 \cr\end{pmatrix}
\end{eqnarray*}
from which it follows that\\ \smallskip \\
\hspace*{.25in}$F_r(x)=\mathlarger{\frac{\partial F}{\partial r}(x)} =
rI - \begin{pmatrix} \hspace{.1in} 0 & \hspace{.1in} 1 & \hspace{.1in} 0 \hspace{.03in} \cr
-1 & \hspace{.1in} 0 & \hspace{.1in} 0\cr
\hspace{.1in} 0 & \hspace{.1in} 0& \hspace{.1in} 0\cr\end{pmatrix}
+ \begin{pmatrix}-2r& \hspace{.2in} 0& \hspace{.1in} t \hspace{.03in}\cr
\hspace{.2in} 0& -2r& \hspace{.1in} s\cr
\hspace{.2in} t& \hspace{.2in} s& \hspace{.1in} 0\cr\end{pmatrix}$\\
\hspace*{.25in}$F_s(x)=\mathlarger{\frac{\partial F}{\partial s}(x)} =
sI - \begin{pmatrix} \hspace{.05in} 0 & \hspace{.1in} 0 & -1 \hspace{.03in} \cr
\hspace{.05in} 0 & \hspace{.1in} 0 & \hspace{.1in} 0\cr
\hspace{.05in} 1 & \hspace{.1in} 0& \hspace{.1in} 0\cr\end{pmatrix}
+ \begin{pmatrix}-2s& \hspace{.1in} t& \hspace{.2in} 0 \hspace{.03in}\cr
\hspace{.2in} t& \hspace{.1in} 0& \hspace{.2in} r\cr
\hspace{.2in} 0& \hspace{.1in} r& -2s\cr\end{pmatrix}$\\
\hspace*{.25in}$F_t(x)=\mathlarger{\frac{\partial F}{\partial t}(x)} =
tI - \begin{pmatrix} \hspace{.05in} 0 & \hspace{.1in} 0 & \hspace{.1in} 0 \hspace{.03in} \cr
\hspace{.05in} 0 & \hspace{.1in} 0 & \hspace{.1in} 1\cr
\hspace{.05in} 0 & -1& \hspace{.1in} 0\cr\end{pmatrix}
+ \begin{pmatrix} \hspace{.1in} 0& \hspace{.2in} s& \hspace{.2in} r \hspace{.03in} \cr
\hspace{.1in} s& -2t& \hspace{.2in} 0\cr
\hspace{.1in} r& \hspace{.2in} 0& -2t\cr\end{pmatrix}.\hspace{.05in}$\\ \medskip \\
With $G(x)=F(x)M-M^TF(x)^T$, for $i,j=1,2,3$, letting $G(x)_{i,j}$ be the entry
of~$G(x)$ in its $i^{th}$ row and $j^{th}$ column, then $g(x)=(G(x)_{12},G(x)_{31},G(x)_{23})$.
Finally, with $G_r(x)$, $G_s(x)$, $G_t(x)$ the partials of $G$ at~$x$, i.e.,
$G_r(x) = F_r(x)M-M^TF_r^T(x)$, $G_s(x) = F_s(x)M-M^TF_s^T(x)$,
$G_t(x) = F_t(x)M-M^TF_t^T(x)$, for $i,j=1,2,3$, letting $G_r(x)_{ij}$, $G_s(x)_{ij}$,
$G_t(x)_{ij}$ be the entries of $G_r(x)$, $G_s(x)$, $G_t(x)$, respectively, in their $i^{th}$
row and $j^{th}$ column, then it is not hard to show that the Jacobian matrix for $g$ at~$x$ is
\begin{eqnarray*}
J(x)=\begin{pmatrix}G_r(x)_{12} & G_s(x)_{12} & G_t(x)_{12}\cr
                    G_r(x)_{31} & G_s(x)_{31} & G_t(x)_{31}\cr
                    G_r(x)_{23} & G_s(x)_{23} & G_t(x)_{23}\cr\end{pmatrix}
\end{eqnarray*}
which is needed for Newton's method.
\par The procedure just described as well as the SVD method carried out as described in the
previous section, have been implemented as part of a Fortran program called {\tt maxtrace.f}, and
this program has been found in our experiments to be close to one hundred percent successful
(it is successful when Newton's method does not fail and therefore the SVD is not used)
on $3\times 3$ nonsymmetric matrices of rank two and three,
but not successful on $3\times 3$ nonsymmetric matrices of rank one. As input to program
{\tt maxtrace.f}, a million $3\times 3$ matrices of random entries were generated and saved in a data
file called randomtrix. With initial point $x_0= (0,0,0)^T$ for each input matrix, program
{\tt maxtrace.f} was then executed on the one million input matrices with an average of 7 to 8
iterations of Newton's method per input matrix that produced solutions for all one million
matrices, i.e., produced one million rotation matrices that transform the one million input
matrices into symmetric matrices. Together with computations also implemented in program
{\tt maxtrace.f} (without the SVD) as described in the previous section,
for obtaining from these symmetric matrices the corresponding one million rotation matrices
that transform them into matrices of maximal trace over rotation matrices, the total time of
the execution of {\tt maxtrace.f} was about $25$~seconds. However, using our Fortran version of the SVD
method only (no Newton's method), Fortran code that is
also part of {\tt maxtrace.f} was also executed that took about $25$~seconds as well for computing
rotation matrices that transform the one million input matrices into the same one million
matrices of maximal trace over rotation matrices obtained with the procedure above. Thus, it
appears that at least for code all written in Fortran, including the SVD method, it takes about
the same amount of time when everything is done using the procedure with Newton's method (and
the SVD method in case Newton's method fails) as it does when everything is done with the SVD
method only. Accordingly, an integer variable called SVDONLY exists in program {\tt maxtrace.f} for
deciding which of the two ways to use. If SVDONLY is set to one, then the latter is used. Otherwise,
if SVDONLY is not set to one, say zero, then the former is used. We also note that in the Fortran code
an integer variable called ITEX exists which is set to the maximum number of allowed iterations
of Newton's method per input matrix.
\,\medskip\\
On the other hand, using Matlab's version of the SVD method only (no Newton's method),
Matlab code under the name {\tt svdcmp.m} was also implemented and executed for computing rotation
matrices that transform the one million input matrices into the same one million matrices of
maximal trace over rotation matrices obtained with the Fortran code above. This was accomplished
in about~$150$ seconds. Actually, program {\tt svdcmp.m}, although a Matlab program, also has the
capability of executing Fortran program {\tt maxtrace.f} to produce the same results obtained above.
This is done with a Matlab mex file called TD$\_$MEX$\_$MAXTRACE.F of {\tt maxtrace.f}. Accordingly,
a Matlab variable called IFLAG exists in program {\tt svdcmp.m} for deciding which to use between Matlab's
SVD method and the Matlab mex file of {\tt maxtrace.f}. If IFLAG is set to one, then the former is used.
Otherwise, if IFLAG is not set to one, say zero, then the latter is used. We note that if IFLAG is not
set to one so that the Matlab mex file of {\tt maxtrace.f} is used, then integer variable SVDONLY described
above must be taken into account as it is part of {\tt maxtrace.f}. Finally we note that with IFLAG equal
to zero, program {\tt svdcmp.m} was successfully executed (using the Matlab mex file of {\tt maxtrace.f})
and took about~$100$ seconds for both SVDONLY equal to zero and equal to one. See Table~1 for a summary
of the times of execution of the Matlab and Fortran codes for the various options
described above. Note the Matlab code is always at least four times slower than the Fortran code.
\,\medskip\\
The Fortran code ({\tt maxtrace.f}), the Matlab code ({\tt svdcmp.m}), the Matlab mex file of
{\tt maxtrace.f} (TD$\_$MEX$\_$MAXTRACE.F), the compiled Matlab mex file of {\tt maxtrace.f}
(TD$\_$MEX$\_$MAXTRACE.mexa64) and a data file consisting of one thousand random $3\times 3$
matrices (randomtrix) can all be obtained at the following links\\
\hspace*{.35in}\verb+https://doi.org/10.18434/M32081+\\
\hspace*{.35in}\verb+http://math.nist.gov/~JBernal+ \verb+/Maximal_Trace.zip+
\begin{table}[H]
        \centering
        \caption{Times of execution in seconds. One million matrices processed.}
        \small
        \begin{tabular}{cc}
                \hline
                Code and options & Time of execution \\ \hline
                Matlab code with Matlab SVD only (no Fortran) & 150 \\
                Matlab code with Fortran mex file and SVD only & 100 \\
                Matlab code with Fortran mex file and Newton's method & 100 \\
                Fortran code with SVD only & 25 \\
                Fortran code with Newton's method & 25 \\
                \hline
        \end{tabular}
%
\end{table}
\,
\section*{\normalsize Summary}
\addcontentsline{toc}{section}{Summary}
In this paper we analyze matrices of maximal trace over rotation matrices. A $d\times d$
matrix~$M$ is of maximal trace over rotation matrices if given any $d\times d$ rotation matrix~$U$,
the trace of $UM$ does not exceed that of~$M$. Given a $d\times d$ matrix~$M$ that is not of maximal
trace over rotation matrices, it is well known that a $d\times d$ rotation matrix~$U$ can be
computed with a method called the Kabsch-Umeyama algorithm (loosely referred to as ``the SVD method''
throughout the paper), based on the computation of the singular value decomposition (SVD) of~$M$
so that $UM$ is of maximal trace over rotation matrices.
Computing a rotation matrix $U$ in this manner for some matrix~$M$ is what is usually
done to solve the constrained orthogonal Procrustes problem and its generalization, Wahba's
problem. As a result of the analysis, we identify a characterization of matrices of maximal
trace over rotation matrices: A $d\times d$ matrix is of maximal trace over rotation matrices
if and only if it is symmetric and has at most one negative eigenvalue, which, if it exists, is
no larger in absolute value than the other eigenvalues of the matrix. Establishing this
characterization is the main goal of this paper, and for $d=2,3$, it is shown how this
characterization can be used to determine whether a matrix is of maximal trace over rotation
matrices. Finally, although depending only slightly on the characterization, as a secondary
goal of the paper, for $d=2,3$, we identify alternative ways, other than the SVD, of
obtaining solutions to the problems of interest.
Given a $2\times 2$ matrix~$M$ that is not of maximal trace over rotation
matrices, an alternative approach that does not involve the SVD method for computing a rotation
matrix~$U$ so that $UM$ is of maximal trace over rotation matrices, is identified that produces
solutions in closed form. Similarly, if $M$ is a $3\times 3$ symmetric matrix, an alternative
approach is also identified that produces solutions partially in closed form. On the other hand,
if $M$ is a $3\times 3$ matrix that is not symmetric, which is the most likely situation when
solving the constrained orthogonal Procrustes problem and Wahba's problem, part of the approach
can still be used to produce the usual orthogonal matrices necessary to carry out the SVD method.
Finally, the situation in which the $3\times 3$ matrix $M$ is not symmetric is reconsidered,
and a procedure is identified that uses the so-called Cayley transform in conjunction with
Newton's method to find a $3\times 3$ rotation matrix $U$ so that $UM$ is symmetric, possibly
of maximal trace over rotation matrices. If the resulting $UM$ is not of maximal trace over
rotation matrices, using the fact that $UM$ is symmetric, another $3\times 3$ rotation matrix
$R$ can then be computed (without the SVD) as described above so that $RUM$ is of maximal trace
over rotation~matrices. Since Newton's method can fail, whenever this happens, as a last resort
the SVD method can then be used also as described above.
We note that all of the above about the three-dimensional case, including the SVD method carried out
as described above, has been successfully implemented in Fortran, and without the SVD, for
randomly generated matrices, the Fortran code is successful in our experiments close to one
hundred percent of the time, using the SVD only when it is not.
Links to the code are provided in the last section of the paper. However, we also note that
it appears that at least for code all written in Fortran, it takes about the same amount of time
when everything is done using the procedure with Newton's method (and the SVD method in case
Newton's method fails) as it does when everything is done with the SVD method only.
We note as well that Matlab code is also provided at the same links, for executing the Fortran
code as a Matlab mex file. Finally, we note that the Matlab code can also be made to compute
solutions using the Matlab version of the SVD method only (no Fortran code executed).
Either way, the Matlab code is at least four times slower than the Fortran code. See Table~1.

\end{document}